\documentclass[a4paper]{elsarticlemod}

\usepackage{amsfonts}
\usepackage{amsmath}
\usepackage{amssymb}
\usepackage{graphicx}
\graphicspath{{figures/}} 
\usepackage{adjustbox}
\usepackage{stmaryrd}
\usepackage{algorithm}
\usepackage{algorithmic}
\usepackage{color}
\usepackage[dvipsnames]{xcolor}
\usepackage[T1]{fontenc}
\usepackage{latexsym}
\usepackage{subfigure}
\usepackage{psfrag}
\usepackage{upgreek}
\usepackage{mathrsfs}
\usepackage{MnSymbol}
\usepackage{float}
\usepackage[colorinlistoftodos, obeyFinal]{todonotes}
\usepackage{textcomp}
\usepackage{soul}
\usepackage{tikz}
\usetikzlibrary{external,shapes.misc}
\usetikzlibrary{positioning,arrows,mindmap,backgrounds, shapes, decorations.fractals, calc, intersections,decorations.pathmorphing,patterns}
\usepackage{pgfplots}
\usepackage{etoolbox}
\usepgfplotslibrary{patchplots}
\usepgfplotslibrary{colormaps}
\usepgfplotslibrary{groupplots}
\usepackage{multirow}
\usepackage{colortbl}
\usepackage{hhline}
\usepackage[title,toc,page]{appendix}

\usepackage{verbatim} 

\def\rmd{{\mathrm{d}}}

\newcommand{\R}{\mathbb{R}}

\newcommand{\eg}{{e.g.,~}}
\newcommand{\CC}{\mathbb{C}}
\newcommand{\CCo}{\overline{\mathbb{C}}}
\newcommand{\DD}{\mathbb{D}}
\newcommand{\DDo}{\overline{\mathbb{D}}}
\newcommand{\BB}{\mathbb{B}}

\newcommand{\T}[1]{{\boldsymbol{#1}}}
\newtheorem{remark}{Remark}

\makeatletter
\newcommand{\vast}{\bBigg@{13.}}
\newcommand{\Vast}{\bBigg@{9}}
\newcommand{\thickhline}{
	\noalign {\ifnum 0=`}\fi \hrule height 1pt
	\futurelet \reserved@a \@xhline
}
\newcolumntype{?}{!{\vrule width 1pt}}
\makeatother
\newcommand\Tstrut{\rule{0pt}{2.6ex}}         
\newcommand\Bstrut{\rule[-0.9ex]{0pt}{0pt}}   

\definecolor{zaffre}{rgb}{0.0, 0.08, 0.66}
\definecolor{bondiblue}{rgb}{0.0,0.58,0.71}
\definecolor{dred}{rgb}{0.92,0,0}
\definecolor{dgreen}{rgb}{0,0.6,0}

\pgfplotsset{stress plot style/.style={
			width=0.48\textwidth,
			y label style={at={(axis description cs:0.15,.5)},anchor=south},
			ylabel={$x_3$ [mm]}}}

\newif\ifrecompiletikz
\recompiletikzfalse
\ifrecompiletikz
  \tikzexternaldisable
\else
  \tikzexternaldisable
\fi

\tikzset{cross/.style={cross out, minimum size=2*(#1-\pgflinewidth), inner sep=0pt, outer sep=0pt},
cross/.default={1pt}}

\newcommand{\blueline}{\textcolor{blue}{$\boldsymbol{\leftrightline}$}}
\newcommand{\redcrosses}{\textcolor{red}{$\boldsymbol{\times}$}}
\DeclareRobustCommand{\blackcircles}{\protect\raisebox{0.0cm}{\tikz{\node[draw=black,line width=0.8pt,scale=0.3,circle,fill=black](){};}}}

\DeclareRobustCommand{\bluesolidcircle}{\protect\raisebox{0.0cm}{\tikz{\node[draw=blue,line width=0.8pt,scale=0.6,circle,fill=none](){}; \draw[draw=blue,line width=0.8pt] (-0.25,0)--(0.25,0);}}}
\DeclareRobustCommand{\blacksolidx}{\protect\raisebox{0.0cm}{\tikz{\node[draw=black,line width=0.8pt,scale=0.6,cross=4pt,fill=none](){}; \draw[draw=black,line width=0.8pt] (-0.25,0)--(0.25,0);}}}
\DeclareRobustCommand{\redsolidcross}{\protect\raisebox{0.0cm}{\tikz{\node[draw=red,line width=0.8pt,scale=0.6,cross=3.5pt,rotate=45,fill=none](){}; \draw[draw=red,line width=0.8pt] (-0.25,0)--(0.25,0);}}}

\begin{document}

\begin{frontmatter}
\title{Accurate equilibrium-based interlaminar stress recovery for isogeometric laminated composite Kirchhoff plates}
\author[pavia]{Alessia Patton\corref{cor1}}
\author[losanne]{Pablo Antol\'in}
\author[arlington]{John-Eric Dufour}
\author[munchen]{Josef Kiendl}
\author[pavia]{Alessandro Reali}
\address[pavia]{Department of Civil Engineering and Architecture - University of Pavia\\
	via Ferrata 3, 27100, Pavia, Italy}  
\address[losanne]{Institute of Mathematics - \'Ecole Polytechnique F\'ed\'erale de Lausanne\\
	CH-1015 Lausanne, Switzerland} 
\address[arlington]{Mechanical and Aerospace Engineering Department - University of Texas at Arlington\\
	500 W 1st St, Arlington, TX 76010} 
\address[munchen]{Department of Civil Engineering and Environmental Sciences - Universit\"at der Bundeswehr M\"unchen\\Werner-Heisenberg-Weg 39, 85577 Neubiberg, Germany}
\cortext[cor1]{Corresponding author. Email: alessia.patton01@universitadipavia.it}

\begin{abstract}

In this paper, we use isogeometric Kirchhoff plates to approximate composite laminates adopting the classical laminate plate theory. Both isogeometric Galerkin and collocation formulations are considered. Within this framework, interlaminar stresses are recovered through an effective post-processing technique based on the direct imposition of equilibrium in strong form, relying on the accuracy and the higher continuity typically granted by isogeometric discretizations. The effectiveness of the proposed approach is proven by extensive numerical tests.

\end{abstract}

\end{frontmatter}

\section{Introduction}\label{sec:introduction}
Laminated composite structures are formed by a collection of laminae stacked to achieve improved mechanical
properties. Each lamina is commonly composed of a matrix that
surrounds and holds in place the fibers, which can be variously oriented giving designers the flexibility
to tailor laminate stiffness and strength still maintaining a reduced weight and matching even demanding structural requirements
(see, e.g., \cite{Gibson1994,Jones1999}). Due to their appealing features, the interest for composite structures in
the engineering community has constantly grown in recent years, especially in the aerospace and automotive industries.

Given the mismatch of material properties of the different layers, laminated composites often exhibit complex behaviors under external loads, which may lead to the typical failure mode referred to as delamination (i.e., separation along layer interfaces). To properly design or assess the structural response of laminated structures, an accurate evaluation of the three-dimensional stress state through the thickness is therefore of paramount importance \cite{Sridharan2008, Mittelstedt2007}. To analyze laminated composite plates, two main categories of approaches are typically identified, namely, two-dimensional \textit{equivalent-single-layer} (ESL) and \textit{layerwise} (LW) theories \cite{Reddy2004,Carrera2014,Liew2019}. Displacement-based ESL theories treat a 3D laminate as an equivalent single-layer plate adopting
suitable kinematics assumptions and therefore implying that displacements are continuous functions of
the thickness coordinate. This results in continuous transverse strains, and, together with ply-wise
discontinuous material properties, necessarily leads to discontinuous through-the-thickness out-of-plane stresses, which violate
what is prescribed by equilibrium. Nevertheless, in addition to inherent simplicity and low computational cost ESL theories can provide a sufficiently accurate description of the global response for thin plates, at least in regions sufficiently far from edges and cut-out boundaries. ESL-based methods include among others the \textit{classical laminate plate theory} (CLPT), which assumes that it is possible to neglect the strains acting through the laminate thickness. Instead, the \textit{first order shear deformation theory} considers the transverse shear strain to be constant with respect to the thickness coordinate and therefore requires shear correction factors, which are difficult to determine for arbitrarily laminated composite plate structures. To overcome this, second- or higher-order ESL laminated plate theories use higher-order polynomials in the expansion of the displacement components through the thickness of the laminate.  However, it is to be noted that higher-order theories introduce additional unknowns that are often difficult to interpret in physical terms and increase the computational effort. On the other hand, full displacement-based LW theories use ply-wise expansions for all three primal variable components such that the 3D displacement field exhibits only $C^0$-continuity through the laminate thickness at the layer interfaces, allowing for a good approximation of interlaminar stresses. The main limitation of  LW theories is that variables are tightly related to the number of layers, leading to high computational costs especially in the case of laminates made of a significant amount of plies. For further approximation theories of laminated composite structures such as the ``Carrera Unified Formulation'' the reader is referred to, e.g., \cite{Carrera2011} and references therein.

Isogeometric analysis (IGA) has been originally introduced in 2005 \cite{Hughes2005} to tightly connect design and analysis, employing shape functions typically belonging to Computer Aided Design (such as B-Splines or NURBS) to approximate both geometry and field variables.
This leads to a cost-saving simplification of the typically expensive mesh generation and refinement processes required by standard finite element analysis. Moreover, the  high smoothness achievable by such functions guarantees superior approximation properties and opens the door to the discretization of high-order PDEs in primal form such as in CLPT, which can be regarded as the extension of Kirchhoff plate theory to laminated composite plates. IGA proved to
be successful in a wide variety of solid and structural problems (see, e.g., the recent works \cite{Leonetti2019,Nitti2020,Antolin2020,Coradello2020} and references therein) and has already been used to solve composite and sandwich plates. 
In particular, IGA has been shown to provide good results when combined with the LW concept (see, e.g., \cite{Guo2014, Guo2015}).
In this context, also 2D isogeometric finite element approaches have been proposed in the literature \cite{Kapoor2013,Nguyen-Xuan2013}, with some of them relying on high-order theories \cite{Thai2015} or employing enhanced shell and plate theories \cite{Remmers2015,Adams2020}.

In this manuscript we present a displacement-based CLPT
approach within the isogeometric analysis framework.
According to this plate theory, interlaminar stresses are identically zero when computed using the constitutive equations. However, these stresses do exist in reality, and they can be responsible for failures in composite laminates because of the difference in the material properties between the layers. Therefore, the proposed modeling strategy is coupled with a post-processing technique which directly relies on equilibrium and grants a highly accurate
prediction of the out-of-plane stress state even from a very coarse 2D displacement solution (e.g., using one high-order element to model the plate mid-plane). 
The adopted post-processing technique takes its origin in \cite{Daghia2008,Engblom1985,Pryor1971} and has already been proven to provide good results for 3D solid plates in the context of both IGA Galerkin \cite{Dufour2018} and collocation  \cite{Patton2019} methods (but also of methods based on Radial Basis Functions \cite{Chiappa2020}).
The effectiveness of the proposed approach relies on the capability to obtain accurate in-plane results with only one element through the thickness and on the possibility to compute accurate stresses and stress derivatives from the obtained displacement field, thanks to the shape function higher-order in-plane continuity properties.

The structure of the paper is as follows. In Section \ref{sec:kirchhoff_plates}, we focus on CLPT basics, considering Kirchhoff plates under bending with ``multiple specially orthotropic layers''. Fundamentals of bivariate B-Splines are presented in Section \ref{sec:numerical_form}, followed by the proposed numerical isogeometric formulations for laminated plates. Such displacement-based modeling strategies do not allow for an immediate assessment of the out-of-plane stress distributions, which can be recovered using an equilibrium-based post-processing technique, as detailed in Section \ref{sec:stress_recovery}. In Section \ref{sec:results}, several numerical tests are considered, showing the ability of the proposed approach to obtain accurate in-plane and out-of-plane stress states. Furthermore we test the behavior of different meshes for increasing length-to-thickness plate ratios and numbers of
layers to show the effectiveness of the method. We also investigate the approach behavior at the plate boundary, where stress concentrations might occur in laminates due to different layerwise material distributions. Finally, conclusions are drawn in Section \ref{sec:conclusions}.

\section{Kirchhoff laminated plates under bending: Layerwise specially orthotropic elasticity}\label{sec:kirchhoff_plates}
In this section we focus on plates with ``multiple specially orthotropic layers'', i.e., laminates characterized by multiple plies for which the bending-stretching coupling coefficients and bending-twisting contributions are zero. This leads the analysis to be greatly simplified because the bending deformation is uncoupled from the extensional deformation \cite{Reddy2004}. Therefore, focusing on the bending case, we acknowledge that the proposed approach is rigorous only for plates characterized by symmetric ply stacking sequences, while for layer arrangements non-symmetric about the mid-plane the coupling phenomenon between bending and stretching is in general not negligible. Neverthless, we will numerically prove in Section \ref{sec:results} that the presented technique is able to provide reasonable approximations to more complex laminates such as antisymmetric cross-ply laminates, namely plates characterized by an even number of layers of equal thickness and the same material properties, with alternating 0\textdegree~and 90\textdegree~orientations.

Under these premises, we recall that the extension of the Kirchhoff plate theory
to laminated composite plates, known as ``classical laminate plate theory'' (CLPT), is based for the bending case on the following displacement field
\begin{subequations}\label{eq:displ_field}
\begin{align}
&u_1(x_1,x_2,x_3)=-x_3w_{,1}\,,\\
&u_2(x_1,x_2,x_3)=-x_3w_{,2}\,,\\
&u_3(x_1,x_2)=w\,,\label{eq:displ_w}
\end{align}
\end{subequations}
where $(u_1,u_2,u_3)$ are the displacement components along the cartesian coordinate directions $(x_1,x_2,x_3)$ of a point belonging to the plate mid-plane (for which $x_3$ is the out-of-plane coordinate) and $w$ is the ``transverse deflection''. The displacement field \eqref{eq:displ_field} implies that straight fibers, normal to the $x_1x_2$-plane before deformation,
remain straight and normal to the mid-surface after deformation.
In equation system \eqref{eq:displ_field} and hereinafter we adopt the convention that the portion of a subscript prior to a comma indicates components, while the portion after the comma refers to partial derivatives; for example, $\sigma_{12,13}=\cfrac{\partial^2 \sigma_{12}}{\partial x_{1}\partial x_{3}}\,$. Small deformations and small strains are assumed throughout the paper.

\subsection{Constitutive relations}\label{subsec:const_rel}
Assuming the displacement field \eqref{eq:displ_field}, the Kirchhoff plate model neglects both transverse shear and membrane strains, while the non-zero corresponding bending strains $\varepsilon_{11}$, $\varepsilon_{22}$, and $\varepsilon_{12}$ cause bending stresses $\sigma_{11}$, $\sigma_{22}$, and $\sigma_{12}$.

Hereinafter, Einstein's notation on repeated indices is used, as well as the convention for which indices in Latin letters take values \{1,2,3\} whereas indices in Greek letters take values \{1,2\}. Accordingly, in-plane strains are defined as 
\begin{equation}\label{eq:strains}
\varepsilon_{\gamma\delta}=-x_3w_{,\gamma\delta}=-x_3\kappa_{\gamma\delta}\,,
\end{equation}
where $\kappa_{\gamma\delta}=w_{,\gamma\delta}$ are the curvatures of the deflected mid-surface and the stress-strain relations for a linear elastic Kirchhoff plate are given by
\begin{equation}\label{eq:stresses}
\sigma_{\alpha\beta}=\CC_{\alpha\beta\gamma\delta}\varepsilon_{\gamma\delta}\,.
\end{equation}

In Section \ref{sec:introduction} we have introduced laminated composite plates as structures made of variously oriented orthotropic elastic plies. For the sake of simplicity, but without loss of generality, we focus here on specially orthotropic layers, for which the principal material coordinates coincide with those of the plate.
Therefore, the number of elastic coefficients of the fourth order elasticity tensor $\mathbb{C}_{ijkl}$ reduces to nine, which, in Voigt notation, can be expressed in terms of engineering constants as
\begin{equation}\label{eq:CIVeltensor}
{\mathbb{C}=\renewcommand\arraystretch{2}\vast{[}\begin{matrix}\mathbb{C}_{11} & \mathbb{C}_{12} & \mathbb{C}_{13} & 0 & 0 & 0\\
& \mathbb{C}_{22} & \mathbb{C}_{23} & 0 & 0 & 0\\
&  & \mathbb{C}_{33} & 0 & 0 & 0\\
&  \text{symm.}&  & \mathbb{C}_{44} & 0 & 0\\
&  &  &  & \mathbb{C}_{55} & 0\\
&  &  &  &  & \mathbb{C}_{66}\\
\end{matrix}\vast{]}=\vast{[}\begin{matrix}
\cfrac{1}{E_{1}} & -\cfrac{\nu_{12}}{E_{1}} & -\cfrac{\nu_{13}}{E_{1}} & 0 & 0 & 0\\
& \cfrac{1}{E_{2}} & -\cfrac{\nu_{23}}{E_{2}} & 0 & 0 & 0\\
&  & \cfrac{1}{E_{3}} & 0 & 0 & 0\\
&  \text{symm.}&  & \cfrac{1}{G_{23}} & 0 & 0\\
&  &  &  & \cfrac{1}{G_{13}} & 0\\
&  &  &  &  & \cfrac{1}{G_{12}}\\
\end{matrix}\vast{]}^{-1}\,.}
\end{equation} We remark that the orthotropic elasticity tensor $\CC$ is not necessarily constant for each ply. Therefore, with $\CC(x_3)$ we denote its through-the-thickness dependency, which is a key aspect in the description of quantities referred to the plate mid-plane.

In accordance with Equation \eqref{eq:stresses} we introduce the bending moments $M_{11}$, $M_{22}$, and $M_{12}$ which are stress resultants with the dimension of moments per unit length
\begin{equation}\label{eq:moments}
M_{\alpha\beta}=\int_{-t/2}^{t/2}x_3\sigma_{\alpha\beta}\rmd x_3\,,
\end{equation}
and, substituting Equation \eqref{eq:strains} into \eqref{eq:stresses}, we combine the obtained
expressions with the bending moment relations \eqref{eq:moments} obtaining
\begin{equation}\label{eq:moments_displ}
M_{\alpha\beta}=-\int_{-t/2}^{t/2}x^2_3\CC_{\alpha\beta\gamma\delta}(x_3)\kappa_{\gamma\delta}\rmd x_3\,,
\end{equation}
where $t$ is the total plate thickness.

Finally, recalling that $\kappa_{\gamma\delta}$ does not depend on the out-of-plane coordinate, we can rewrite \eqref{eq:moments_displ} as
\begin{equation}\label{eq:moments_displ_compact}
M_{\alpha\beta}=-\DD_{\alpha\beta\gamma\delta}\kappa_{\gamma\delta}\,,
\end{equation}
being $\DD_{\alpha\beta\gamma\delta}$ the bending material stiffness, defined as
\begin{equation}\label{eq:DD}
\DD_{\alpha\beta\gamma\delta}=\int_{-t/2}^{t/2}x^2_3\CC_{\alpha\beta\gamma\delta}(x_3)\rmd x_3\,.
\end{equation}

\subsection{Boundary-value problem}\label{subsec:bvp}
The boundary value problem associated with an elastic Kirchhoff plate under bending can be formulated as follows.

Let $\Omega$ be an open subset of $\mathbb{R}^2$, subjected to a transversal, i.e., normal to the plate mid-plane, load $q:\Omega\mapsto\R$. We assume that $\Omega$ has a sufficiently smooth boundary $\Gamma$ with a well-defined normal $\boldsymbol{n}$. $\Gamma$ can be decomposed as $\Gamma=\overline{\Gamma_w\cup\Gamma_Q}$ and $\Gamma=\overline{\Gamma_\varphi\cup\Gamma_M}$ with $\Gamma_w\not=\emptyset$ and $\Gamma_w\cap\Gamma_Q=\emptyset$, $\Gamma_\varphi\cap\Gamma_M=\emptyset$. Given the distributed load $q$, and the boundary condition functions $w_\Gamma:\Gamma_w\mapsto\mathbb{R}$, $\varphi_\Gamma:\Gamma_\varphi\mapsto\mathbb{R}$, $Q_\Gamma:\Gamma_Q\mapsto\mathbb{R}$, $M_\Gamma:\Gamma_M\mapsto\mathbb{R}$, we look for the transverse deflection $w:\Omega\mapsto\mathbb{R}$ such that
\begin{subequations}\label{eq:plate_bvp}
\begin{alignat}{3}
& M_{\alpha\beta,\alpha\beta}=q                                   \quad\quad &&    \quad && \quad\hbox{in}\quad\Omega
\label{eq:plate_pde}\\
& M_{\alpha\beta}=-\DD_{\alpha\beta\gamma\delta}(x_3)\kappa_{\gamma\delta}=-\DD_{\alpha\beta\gamma\delta}(x_3)w_{,\gamma\delta}   \quad\quad &&    \quad && \quad\hbox{in}\quad\Omega\label{eq:plate_constlaw}\\
& w = w_\Gamma                                                    \quad\quad &&    \quad && \quad\hbox{on}\quad\Gamma_w\label{eq:plate_dbc}\\
& w_{,\alpha}n_{\alpha}=\varphi_\Gamma                             \quad\quad &&    \quad && \quad\hbox{on}\quad\Gamma_\varphi\label{eq:plate_rbc}\\
& (M_{\alpha\beta,\beta} + M_{\alpha\delta,\delta})n_{\alpha}=Q_\Gamma \quad\quad && \text{with}\quad\delta\neq\alpha\quad && \quad\hbox{on}\quad\Gamma_Q\label{eq:plate_sbc}\\
& n_{\alpha}M_{\alpha\beta}n_{\beta} = M_\Gamma \quad\quad && \quad && \quad\hbox{on}\quad\Gamma_M\,,\label{eq:plate_mbc}
\end{alignat}
\end{subequations} 
where $w_\Gamma$, $\varphi_\Gamma$, and $M_\Gamma$ represent, respectively, the prescribed normal out-of-plane displacement, rotation, and moment. 
Instead, $Q_\Gamma$ stands for the normal component of the so-called ``effective shear'' (see \cite{Felippa2017,Reali2015,Wang2000}), classically defined by the combination of the effect on the boundary of shear forces
(i.e., $M_{\alpha\beta,\beta}n_{\alpha}$) and twisting moments (i.e., $M_{\alpha\delta,\delta}n_{\alpha}$).

\subsection{Weak form}\label{subsec:weak_form}
In a variational approach, the governing equations are obtained by the principle of virtual displacements.

A given mechanical system can take many possible configurations in accordance with its geometric constraints. Of all the admissible configurations (i.e., the set of configurations that satisfy the geometric constraints), only one also satisfies equilibrium. These configurations can be regarded as infinitesimal variations, during which the compatibility constraints of the system are not violated. Such variations are called virtual displacements and do not have any relation to the actual displacements that might occur due to a change in the applied loads \cite{Reddy2004}. Thus, for a plate occupying a region $\Omega\subseteq\mathbb{R}^2$ and subjected to pure bending, the only contribution to the internal virtual work (in the primal field virtual transverse displacement, $\delta w$) is given by the in-plane bending moments, $M_{\alpha\beta}$, and their relative virtual work conjugate curvatures, $\delta\kappa_{\alpha\beta}$, as
\begin{equation}\label{eq:W_int}
\delta W_{int}[\delta w]=\int_{\Omega}M_{\alpha\beta}\delta\kappa_{\alpha\beta}\rmd\Omega\,.
\end{equation}
The external virtual work is given
instead by the sum of three components \cite{Felippa2017}. These are respectively due to applied lateral loads, $\delta W_{ext,q}$, applied edge moments and
transverse shears, $\delta W_{ext,B}$, and to corner loads, $\delta W_{ext,C}$, i.e.,
\begin{equation}\label{eq:W_ext}
\delta W_{ext}[\delta w]= \delta W_{ext,q}[\delta w] + \delta W_{ext,B}[\delta w] + \delta W_{ext,C}[\delta w]\,.
\end{equation}
The first two terms read
\begin{equation}\label{eq:W_extq}
\delta W_{ext,q}[\delta w]=\int_{\Omega} q\delta w \rmd\Omega\,,
\end{equation}
and
\begin{equation}\label{eq:W_extB}
\delta W_{ext,B}[\delta w]=\int_{\Gamma}(Q_{\Gamma}\delta w_{\Gamma} + M_{\Gamma}\delta\phi_{\Gamma})\rmd\Gamma\,.
\end{equation}
Finally, if the plate has $n_c$ corners at which the displacement $w_j$, with $j=1,2,...,n_c$, is not prescribed, the term $\delta W_{ext,C}$ comes into play considering the so called ``corner forces'', i.e., jumps in the corresponding twisting moments.
In this work we assume for the sake of simplicity that the transverse displacement of a corner is always prescribed, which grants that the contribution of that corner to the external virtual work vanishes because its corresponding displacement variation is zero. This assumption does not constitute any limitation to the purpose of the present work.

\section{Numerical formulations}\label{sec:numerical_form}
In this section we introduce the notions of bivariate B-Splines and detail the numerical isogeometric formulations to approximate the problem variables and thus the equations governing the laminated Kirchhoff plate.

\subsection{Bivariate B-Splines}\label{subsec:2Bsplines}
We introduce the basic definitions and notations regarding bivariate B-Splines, while, for further details, readers may refer to \cite{Cottrell2007,Hughes2005,Piegl1997} and
references therein.

To this end, we need to first define two univariate knot vectors, i.e., non-decreasing set of coordinates in the $d$-th parameter
space, as
\begin{alignat}{2}\label{eq:knotvectors}
&\Xi^{d} = \{\xi_{1}^{d},...,\xi_{m_{d}+p_{d}+1}^{d}\} \quad &&\quad d = 1,2\,,
\end{alignat}
where $p_{d}$ represents the polynomial degree in the parametric direction $d$,
and $m_{d}$ is the associated number of basis functions. A univariate B-Spline basis function $N^{d}_{i_{d},p_{d}}(\xi^{d})$, corresponding to the parametric coordinate $\xi^{d}$, can be then constructed, for each $i_d$ position in the tensor
product structure, using the Cox-de Boor formula starting from $p_{d}=0$
\begin{alignat}{2}\label{eq:univarBSplinesConst}
&N^{d}_{i_{d},p_{d}}(\xi^{d})=\begin{cases}1\quad &\quad\xi_{i_{d}}^{d}\le\xi^{d}<\xi_{i_{d}+1}^{d}\\
0\quad&\quad \text{otherwise}
\end{cases}\,,
\end{alignat}
while the basis functions for $p_{d}>0$ are recursively obtained as
\begin{equation}\label{eq:univarBSplines}
N^{d}_{i_{d},p_{d}}(\xi^{d})=\cfrac{\xi^{d}-\xi^{d}_{i_{d}}}{\xi^{d}_{i_{d}+p_{d}}-\xi^{d}_{i_{d}}}N^{d}_{i_{d},p_{d}-1}(\xi^{d})+\cfrac{\xi^{d}_{i_{d}+p_d+1}-\xi^{d}}{\xi^{d}_{i_{d}+p_d+1}-\xi^{d}_{i_{d}+1}}N^{d}_{i_{{d}+1},p_{d}}(\xi^{d})\,,
\end{equation}
where the convention $0/0=0$ is assumed.

Bivariate basis functions $B_{\textbf{i},\textbf{p}}(\boldsymbol{\xi})$ are obtained by tensor product of two sets of univariate B-Splines as
\begin{equation}\label{eq:multivarBSplines}
B_{\textbf{i},\textbf{p}}(\boldsymbol{\xi})=\prod\limits_{d=1}^{2}N_{i_{d},p_{d}}(\xi^{d})\,,
\end{equation}
where $\textbf{i} = \{i_{1}, i_{2}\}$ plays the role of a multi-index which describes the considered position in the tensor product structure, $\textbf{p} = \{{p_{1},p_{2}}\}$ indicates the polynomial degrees, and $\boldsymbol{\xi} = \{\xi^{1},\xi^{2}\}$ represents the vector of the parametric coordinates in each parametric direction $d$.

Finally, B-Spline bidimensional geometries are built as a linear combination of bivariate B-Spline basis functions as follows
\begin{equation}\label{eq:multivarBSplinesGeom}
\textbf{S}(\boldsymbol{\xi})=\sum\limits_{\textbf{i}}B_{\textbf{i},\textbf{p}}(\boldsymbol{\xi})\textbf{P}_\textbf{i}\,,
\end{equation}
where the coefficients $\textbf{P}_\textbf{i}\in\mathbb{R}^{2}$ are the so-called
control points, and the summation is extended to all combinations of the multi-index \textbf{i}.

\subsection{Constitutive relations: Laminated composite material}
To capture the laminated composite through-the-thickness behavior, we need to account for the proper material distribution layer by layer even though the Kirchhoff theory assumes that a mid-surface plane can be used to represent the three-dimensional solid plate in a two-dimensional form. In order to include the complete ply stacking sequence contribution, we consider the needed 3D material tensor \eqref{eq:CIVeltensor} components for each $k$-th layer and, to create an equivalent single bivariate plate, we homogenize the material properties according to \cite{Sun1988} by means of the following relations
\begin{subequations}\label{eq:aveC}
	\begin{alignat}{2}
	&\CCo_{ab}=\sum_{k=1}^{N}\overline{t}_{k}\CC_{ab}^{(k)}+\sum_{k=2}^{N}(\CC_{a3}^{(k)}-\CCo_{a3})\overline{t}_{k}\frac{(\CC_{b3}^{(1)}-\CC_{b3}^{(k)})}{\CC_{33}^{(k)}}\quad&&\quad a,b=1,2\,,\label{eq:aveC1}\\
	&\CCo_{66}=\sum_{k=1}^{N}\overline{t}_{k}\CC_{66}^{(k)}\,,\quad&&\quad\label{eq:aveC9}
	\end{alignat}
\end{subequations}
where $\overline{t}_{k}=\cfrac{t_{k}}{t}$ represents the volume fraction of the $k$-th lamina, $t$ being the total plate thickness, and $t_{k}$ the $k$-th ply thickness.

At this point $\CCo$ \footnote{We note that in order to obtain $\CCo$ the out-of-plane shear moduli are not considered in accordance with the homogenization rule in \cite{Sun1988} adapted for a bivariate case.} is independent of the $x_3$ coordinate and we can recover the homogenized bending material stiffness \footnote{We would like to underline that from here on out all the presented numerical strategies and results refer to the obtained homogenized bending material stiffness $\DDo$.} from Equation \eqref{eq:DD} as
\begin{equation}\label{eq:DDbar}
\DDo=\begin{pmatrix}\DDo_{11}&\DDo_{12}&0\\
&\DDo_{22}&0\\
\text{symm.}&&\DDo_{66}\end{pmatrix}=\cfrac{t^3}{12}\begin{pmatrix}\CCo_{11}&\CCo_{12}&0\\
&\CCo_{22}&0\\
\text{symm.}&&\CCo_{66}\end{pmatrix}\,.
\end{equation}

\subsection{Isogeometric collocation method}\label{subsec:collocation}
The collocation method can be seen as a Petrov-Galerkin method where the test functions are smoothed Dirac delta functions (converging to the Dirac delta distributions located at the collocation points as the smoothing parameter tends to zero).
It can be therefore regarded as a sort of stable one-point quadrature Galerkin method giving raise to a strong-form method. As reported in \cite{Auricchio2010b}, a delicate issue for collocation methods is the determination of suitable collocation points. The simplest and most widespread approach is to collocate the governing strong-form equations at the images of Greville abscissae (see, \eg\cite{Johnson2005}) and this is the strategy also herein adopted. Accordingly, along each parametric direction $d=1,2$, we consider a set of $m_d$ Greville abscissae, i.e., points obtained from the knot vector components, $\theta^d_{i}$, as 
\begin{alignat}{2}\label{eq:greville}
&\overline{\theta}^d_{i}=\frac{\theta^d_{i+1}+\theta^d_{i+2}+...+\theta^d_{i+p_d}}{p_d}\quad&&\quad i = 1,...,m_d\,,
\end{alignat}
$p_d$ being the degree of approximation and $m_d$ the number of basis functions.
Having defined $\boldsymbol{\tau}$ as the collocation points matrix, such that each $ij$-th entry is $\boldsymbol{\tau}_{ij}=\bigg(\cfrac{\sum_{k=1}^{p_1}\xi_{i+k}}{p_1},\cfrac{\sum_{l=1}^{p_2}\eta_{j+l}}{p_2}\bigg)$ with $i = 1,...,m_1,\;j = 1,...,m_2$, we approximate the displacement field $\textbf{w}$ as a linear combination of bivariate shape functions and control variables $\hat{\textbf{w}}_{\textbf{i}}$ as 
\begin{equation}\label{eq:displapprox}
\textbf{w}(\boldsymbol{\tau})=B_{\textbf{i},\textbf{p}}(\boldsymbol{\tau})\hat{\textbf{w}}_{\textbf{i}}\,.
\end{equation} Following \cite{Reali2015}, without loss of generality, we describe our collocation strategy for the case of a simply supported plate, that is, $\Gamma_w=\Gamma_M=\Gamma$.

In Voigt notation, we can rewrite Equation \eqref{eq:moments_displ_compact} as
\begin{equation}\label{eq:momentsCompact}
\textbf{M}=-\DDo\boldsymbol{\kappa}\,,
\end{equation}
where the bending moment vector is equal to
\begin{equation}\label{eq:bending_moments_vec}
\textbf{M}=\begin{bmatrix} M_{11} & M_{22} & M_{12} \end{bmatrix}^{T}
\end{equation}
and the curvature vector $\boldsymbol{\kappa}$ is defined as
\begin{equation}\label{eq:bending_curvatures_vec}
\boldsymbol{\kappa}=\begin{bmatrix} \kappa_{11} & \kappa_{22} & 2\kappa_{12} \end{bmatrix}^{T}=\begin{bmatrix} w_{,11} & w_{,22} & 2w_{,12} \end{bmatrix}^{T}\,.
\end{equation}
We then insert the approximate displacements~\eqref{eq:displapprox} into the bending moment equations~\eqref{eq:momentsCompact} and we further substitute into equilibrium equations \eqref{eq:plate_pde}, obtaining
\begin{alignat}{2}\label{eq:equilibrium_collocated}
&-\textbf{K}(\boldsymbol{\tau})\hat{\textbf{w}}_{\textbf{i}}=\textbf{q}(\boldsymbol{\tau})\quad&&\quad\forall\boldsymbol{\tau}_{ij}\in\Omega\,,
\end{alignat}
where $\textbf{K}(\boldsymbol{\tau})$ can be expressed as
\begin{equation}\label{eq:Kcoeffs}
\textbf{K}(\boldsymbol{\tau})=\DDo_{11}\cfrac{\partial^4{B_{\textbf{i},\textbf{p}}(\boldsymbol{\tau})}}{\partial{x_1}^4}+
2\bigg(\DDo_{12}+2\DDo_{66}\bigg)\cfrac{\partial^4{B_{\textbf{i},\textbf{p}}(\boldsymbol{\tau})}}{\partial{x_1}^2\partial{x_2}^2}+\DDo_{22}\cfrac{\partial^4{B_{\textbf{i},\textbf{p}}(\boldsymbol{\tau})}}{\partial{x_2}^4}\,,
\end{equation}
while, substituting in \eqref{eq:plate_mbc}, we obtain instead
\begin{alignat}{2}\label{eq:boundaryM_collocated}
&-\tilde{\textbf{K}}(\boldsymbol{\tau})\hat{\textbf{w}}_{\textbf{i}}=\textbf{M}_\Gamma(\boldsymbol{\tau})\quad&&\quad\forall\boldsymbol{\tau}_{ij}\in\Gamma_{M}\,,
\end{alignat}
with $\tilde{\textbf{K}}(\boldsymbol{\tau})$ having the following form
\begin{equation}
\begin{aligned}\label{eq:Ktildacoefs}
\tilde{\textbf{K}}(\boldsymbol{\tau})&=\DDo_{11}\cfrac{\partial^2{B_{\textbf{i},\textbf{p}}(\boldsymbol{\tau})}}{\partial{x_1}^2}n^2_1+
\DDo_{12}\bigg(\cfrac{\partial^2{B_{\textbf{i},\textbf{p}}(\boldsymbol{\tau})}}{\partial{x_2}^2}n^2_1+\cfrac{\partial^2{B_{\textbf{i},\textbf{p}}(\boldsymbol{\tau})}}{\partial{x_1}^2}n^2_2\bigg)\\
&+\DDo_{22}\cfrac{\partial^2{B_{\textbf{i},\textbf{p}}(\boldsymbol{\tau})}}{\partial{x_2}^2}n^2_2+4\DDo_{66}\cfrac{\partial^2{B_{\textbf{i},\textbf{p}}(\boldsymbol{\tau})}}{\partial{x_1}\partial{x_2}}n_1n_2\,.
\end{aligned}
\end{equation}
Regarding boundary condition imposition, the strategy is exactly the same as thoroughly discussed by Reali and Gomez for an isotropic plate, and we therefore refer interested readers to \cite{Reali2015} for futher details.

\subsection{Isogeometric Galerkin method}\label{subsec:galerkin}
For an isogeometric Galerkin approach the variation of the energy functional in a system can be regarded as the sum of all its element-wise variations, thus
\begin{alignat}{1}
&\delta W_{int}=\sum_{e}^{N_{e}}\delta W^{(e)}_{int}\,,\label{eq:Welint}\\
&\delta W_{ext}=\sum_{e}^{N_{e}}\delta W^{(e)}_{ext}\,,\label{eq:Welext}
\end{alignat}
where $N_{e}$ denotes the number of elements in the plate domain and the superscript $(e)$ is the
element index.
Then, approximating the displacement field as a linear combination of bivariate shape functions and control variables as 
\begin{subequations}\label{eq:displapproxGal}
\begin{align}
&\boldsymbol{w}^{(e)}(\boldsymbol{\bar{\xi}})=B^{(e)}_{\textbf{i},\textbf{p}}(\boldsymbol{\bar{\xi}})\hat{\boldsymbol{w}}^{(e)}_{\textbf{i}}\,,\label{eq:displapproxGalw}\\
&\delta\boldsymbol{w}^{(e)}(\boldsymbol{\bar{\xi}})=B^{(e)}_{\textbf{i},\textbf{p}}(\boldsymbol{\bar{\xi}})\delta\hat{\boldsymbol{w}}^{(e)}_{\textbf{i}}\,,\label{eq:displapproxGaldeltaw}
\end{align}
\end{subequations}
we substitute \eqref{eq:displapproxGal} into \eqref{eq:Welint} obtaining the aproximate element internal energy variation 
\begin{equation}\label{eq:W_int_discr}
\delta W^{(e)}_{int}=(\delta \hat{\boldsymbol{w}}^{(e)}_{\textbf{i}})^{T}\int_{\Omega^{(e)}}\BB^{(e)T}_{\textbf{i},\textbf{p}}\,\DDo\,\BB^{(e)}_{\textbf{j},\textbf{p}}\rmd\Omega^{(e)}\hat{\boldsymbol{w}}^{(e)}_{\textbf{j}}\simeq(\delta \hat{\boldsymbol{w}}^{(e)}_{\textbf{i}})^{T}\mathbb{K}^{(e)}\hat{\boldsymbol{w}}^{(e)}_{\textbf{j}}\,,
\end{equation}
where $\mathbb{K}^{(e)}$ is the stiffness matrix computed approximating the integral with a quadrature rule. In this work we consider standard Gauss integration. In Equation \eqref{eq:W_int_discr} $\BB^{(e)}_{\textbf{i},\textbf{p}}$ is defined as
\begin{equation}\label{eq:Bmatrix}
\BB^{(e)}_{\textbf{i},\textbf{p}}=\begin{bmatrix}\cfrac{\partial^2{B^{(e)}_{\textbf{i},\textbf{p}}(\boldsymbol{\bar{\xi}})}}{\partial{x_1}^2}\\
\cfrac{\partial^2{B^{(e)}_{\textbf{i},\textbf{p}}(\boldsymbol{\bar{\xi}})}}{\partial{x_2}^2}\\
2\cfrac{\partial^2{B^{(e)}_{\textbf{i},\textbf{p}}(\boldsymbol{\bar{\xi}})}}{\partial{x_1}\partial{x_2}}\end{bmatrix}\,,
\end{equation}
where $\boldsymbol{\bar{\xi}}$ is the matrix of the quadrature point positions.

\section{Stress recovery procedure}\label{sec:stress_recovery}
Since Kirchhoff theory is intrinsically two-dimensional, assessments of out-of-plane stress distributions are not immediately possible. Thus, the strategies proposed in Sections \ref{subsec:collocation} and \ref{subsec:galerkin} are by themselves not suitable for the calculation of interlaminar stresses but can be easily coupled with an a-posteriori step based on equilibrium, which, following \cite{Dufour2018,Patton2019}, has already proved to allow a rigorous layerwise reconstruction of out-of-plane stresses for laminated solid plates in the context of both isogeometric Galerkin and collocation methods.

The starting point is the fact that stresses must satisfy the equilibrium equations
\begin{equation}\label{eq:equilibrium_compact}
\T{\nabla} \cdot \T{\sigma} + \T{b} = \boldsymbol{0}
\end{equation}
at every point, where $\T{\nabla}\cdot$ represents the divergence operator. Equilibrium equations
\eqref{eq:equilibrium_compact} can be expressed in a componentwise way as:
\begin{subequations}\label{eq:equilibriumeng}
	\begin{align}
	&\sigma_{11,1}+\sigma_{12,2}+\sigma_{13,3}=-b_{1}\,,\label{eq:equilibriumeng1}\\
	&\sigma_{12,1}+\sigma_{22,2}+\sigma_{23,3}=-b_{2}\,,\label{eq:equilibriumeng2}\\
	&\sigma_{13,1}+\sigma_{23,2}+\sigma_{33,3}=-b_{3}\,,\label{eq:equilibriumeng3}
	\end{align}
\end{subequations}
and integrating Equations~\eqref{eq:equilibriumeng1} and~\eqref{eq:equilibriumeng2} along the thickness, we can recover the out-of-plane shear stresses as
\begin{subequations}\label{eq:recoveryShear}
	\begin{align}
	\label{eq:ppsigma13}
	\sigma_{13}(x_3) &= -\int^{x_3}_{\bar{x}_3}(\sigma_{11,1}(\zeta) + \sigma_{12,2}(\zeta)+b_1(\zeta))\rmd \zeta + \sigma_{13}(\bar{x}_3)\,,\\
	\label{eq:ppsigma23}
	\sigma_{23}(x_3) &= -\int^{x_3}_{\bar{x}_3}(\sigma_{12,1}(\zeta) + \sigma_{22,2}(\zeta)+b_2(\zeta))\rmd \zeta + \sigma_{23}(\bar{x}_3)\,,
	\end{align}
\end{subequations}
where $\zeta$ represents the coordinate along the plate thickness direction. Note that in this work all integrals along this direction are computed using a composite trapezoidal quadrature rule.

Once we substitute the appropriate derivatives of the out-of-plane shear stresses \eqref{eq:recoveryShear} into Equation \eqref{eq:equilibriumeng3}, integrating along the thickness, we can recover also $\sigma_{33}$ as:
\begin{equation}\label{eq:recoveryNormal}
\begin{aligned}
\sigma_{33}(x_3) & = \int^{x_3}_{\bar{x}_3}\left[\int^{\zeta}_{\bar{x}_3}(\sigma_{11,11}(\xi) + \sigma_{22,22}(\xi) + 2\sigma_{12,12}(\xi) + b_{1,1}(\xi)+b_{2,2}(\xi))\rmd \xi\right]\rmd \zeta \\ 
& - \int^{x_3}_{\bar{x}_3} b_3(\zeta)  \rmd \zeta - (x_3-\bar{x}_3)(\sigma_{13,1}(\bar{x}_3) + \sigma_{23,2}(\bar{x}_3)) + \sigma_{33}(\bar{x}_3)\,,
\end{aligned}
\end{equation} 
where the integral constants should be chosen to fulfill the boundary conditions at the top or bottom surfaces $\bar{x}_3$ \cite{Dufour2018}.

The derivatives of the in-plane stress components necessary for the proposed post-processing step are computed, from an even very coarse displacement solution as
\begin{subequations}\label{eq:der_der2_sigma}
\begin{align}	
&\sigma_{\alpha\beta,\gamma}=\CC_{\alpha\beta \zeta\eta}(x_3)(-x_3\kappa_{\zeta\eta,\gamma})=\CC_{\alpha\beta \zeta\eta}(x_3)(-x_3w_{,\zeta\eta\gamma})\,,\label{eq:dersigma}\\
&\sigma_{\alpha\beta,\gamma\delta}=\CC_{\alpha\beta \zeta\eta}(x_3)(-x_3\kappa_{\zeta\eta,\gamma\delta})=\CC_{\alpha\beta \zeta\eta}(x_3)(-x_3w_{,\zeta\eta\gamma\delta})\,.\label{eq:der2sigma}
\end{align}
\end{subequations}
From Equations \eqref{eq:dersigma} and \eqref{eq:der2sigma} it is clear that the proposed post-processing step requires the shape functions to be highly continuous (i.e., $C^4$-continuous), which can be easily achieved by means of isogeometric analysis. In the following section, we show convincing numerical experiments proving that such a post-processing technique works nicely in the context of both isogeometric Galerkin and collocation methods.

\section{Numerical Results}\label{sec:results}
In this section, several examples are presented for the static analysis of composite laminated Kirchhoff plates under bending. To this extent we validate the considered tests against Pagano's analytical solution \cite{Pagano1970} to showcase the accuracy of the
proposed post-processing technique in reconstructing the out-of-plane stress field, addressing different aspects such as the method sensitivity to parameters of interest (i.e., number of layers and length-to-thickness
ratio).

\subsection{Analytical solution: The Pagano test case}
\begin{figure}[!htbp]
	\centering
	\includegraphics[width=.5\textwidth]{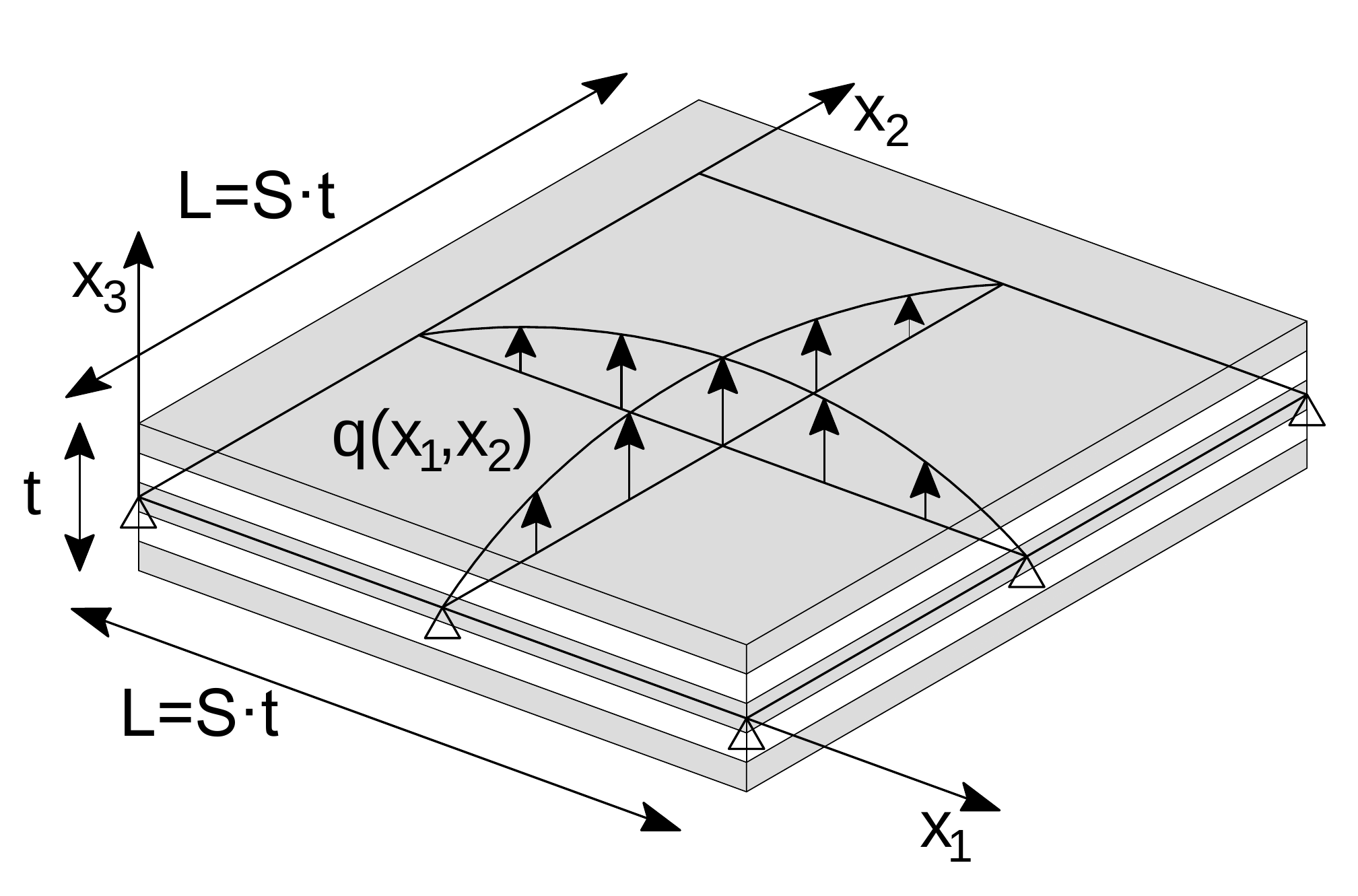}
	\caption{The Pagano test case~\cite{Pagano1970}. Problem geometry.}
	\label{fig:testproblem}
\end{figure}

The Pagano test case considers a solid cross-ply plate of total thickness $t$, made of $N$ orthotropic layers. The structure is simply supported on all edges and subjected to a transverse sinusoidal loading $q(x_1,x_2)$, on the top surface, while the bottom one is traction-free. The thickness of every single layer is set to 1 mm, and the edge length, $L$, is chosen to be $S$ times larger than the total thickness of the laminate.  
We approximate Pagano's solid benchmark with a bivariate plate as in Figure \ref{fig:testproblem}.
\newline
The layer material parameters taken into account for all numerical tests are summarized in Table~\ref{tab:matproperties} for 0\textdegree-oriented plies, while the considered loading pressure, applied at the plate mid-plane, is equal to
\begin{equation}\label{eq:load}
q(x_1,x_2)=\sigma_0\sin(\dfrac{\pi x_1}{St})\sin(\dfrac{\pi x_2}{St})\,,
\end{equation}
where $\sigma_0 =$ 1 MPa.
\begin{table}[!htbp]
	\caption{Adopted material properties for 0\textdegree-oriented layers.} 
	\begin{center}
		\small{
			\begin{tabular}{?c c c c c c c c c?}
				\thickhline 
				$E_{1}$ & $E_{2}$& $E_{3}$ & $G_{23}$ & $G_{13}$ & $G_{12}$ & $\nu_{23}$ & $\nu_{13}$ & $\nu_{12}$\Tstrut\Bstrut\\\hline 
				[GPa] & [GPa] & [GPa] & [GPa] & [GPa] & [GPa]& [-] & [-] & [-]\Tstrut\Bstrut\\\thickhline 
				25000 & 1000 & 1000 & 200 & 500 & 500 & 0.25 & 0.25& 0.25\Tstrut\Bstrut\\\thickhline  
		\end{tabular}}
	\end{center}\label{tab:matproperties}
\end{table}
\newline
With reference to Equations \eqref{eq:plate_mbc} and \eqref{eq:plate_dbc}, the simply supported edge conditions are taken as
\begin{alignat}{2}\label{eq:DBCs}
&M_{\Gamma}=0\quad\hbox{and}\quad w_{\Gamma}=0\quad&&\hbox{on}\quad \Gamma_w=\Gamma_M=\Gamma\,.
\end{alignat}
We remark that for collocation the boundary condition $M_{\Gamma}=0$ is strongly imposed, while it is naturally satisfied in Galerkin methods.

Finally, all results hereinafter reported are expressed in terms of normalized stress components as
\begin{subequations}
\begin{alignat}{2}\label{eq:normalizedresults}
&\bar{\sigma}_{ij}=\dfrac{\sigma_{ij}}{\sigma_0S^2}\quad&&\quad i,j=1,2\,,\\
&\bar{\sigma}_{i3}=\dfrac{\sigma_{i3}}{\sigma_0S}\quad&&\quad i=1,2\,,\\
&\bar{\sigma}_{33}=\dfrac{\sigma_{33}}{\sigma_0}\,.&&
\end{alignat}
\end{subequations}

\subsection{In-plane solution assessment and out-of-plane reconstruction from equilibrium}
In this section, we present and comment several numerical examples considering a cross-ply distribution of layers, namely a 90\textdegree/0\textdegree~stacking sequence from the bottom to the top of the plate. All numerical simulations are carried out using an in-plane degree of approximation $p=q=6$, which fulfills the continuity requirements described in Section \ref{sec:stress_recovery}, and a very coarse grid comprising of 7x7 control points, or equivalently degrees of freedom, which corresponds to only one element (which has been verified to grant good results for this problem).

\begin{figure}[!htbp]
	\centering
	\subfigure[Normalized $\sigma_{11}$\label{subfig-1:in_plane_sig_1_1_l11}]{\ifrecompiletikz\tikzsetnextfilename{fig_02_a}\tikzexternalenable\input{images/fig_02_a}\tikzexternaldisable\else\includegraphics{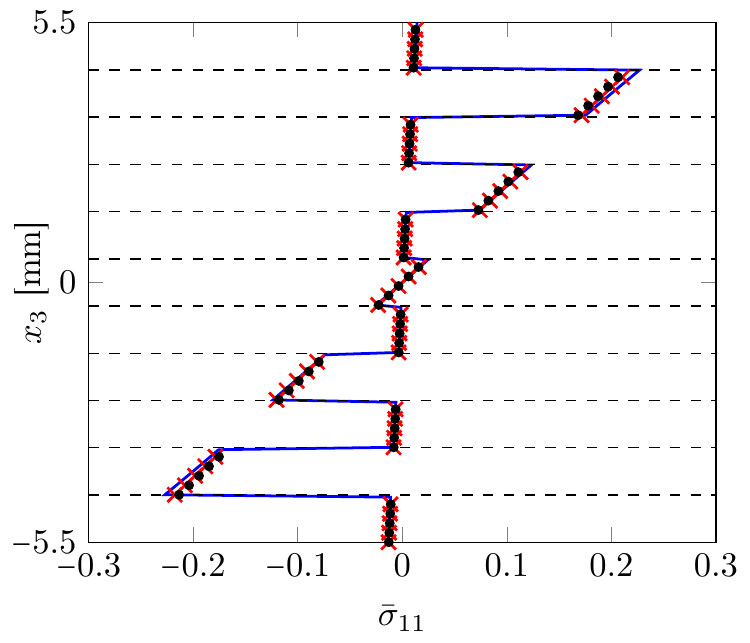}\fi}\hfill
	\subfigure[Normalized $\sigma_{11}$\label{subfig-2:in_plane_sig_1_1_l34}]{\ifrecompiletikz\tikzsetnextfilename{fig_03_a}\tikzexternalenable\input{images/fig_03_a}\tikzexternaldisable\else\includegraphics{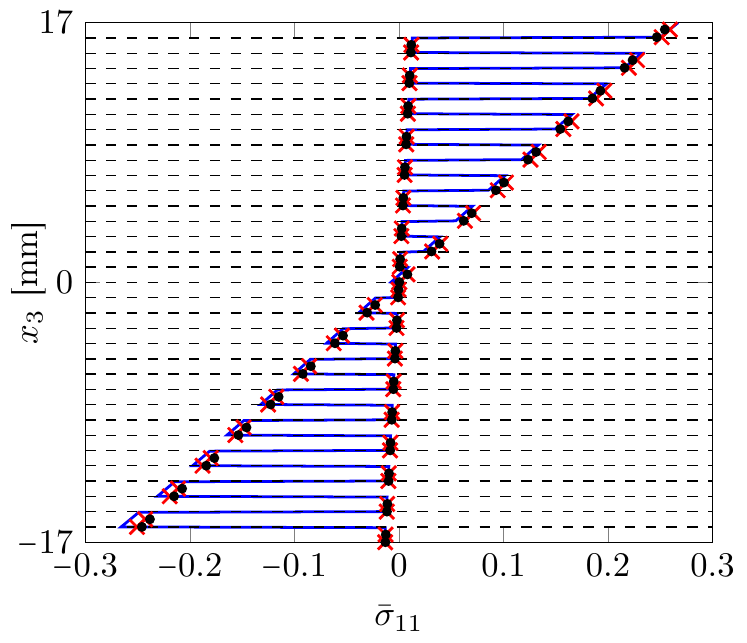}\fi}\\
	\subfigure[Normalized $\sigma_{22}$\label{subfig-3:in_plane_sig_1_1_l11}]{\ifrecompiletikz\tikzsetnextfilename{fig_02_b}\tikzexternalenable\input{images/fig_02_b}\tikzexternaldisable\else\includegraphics{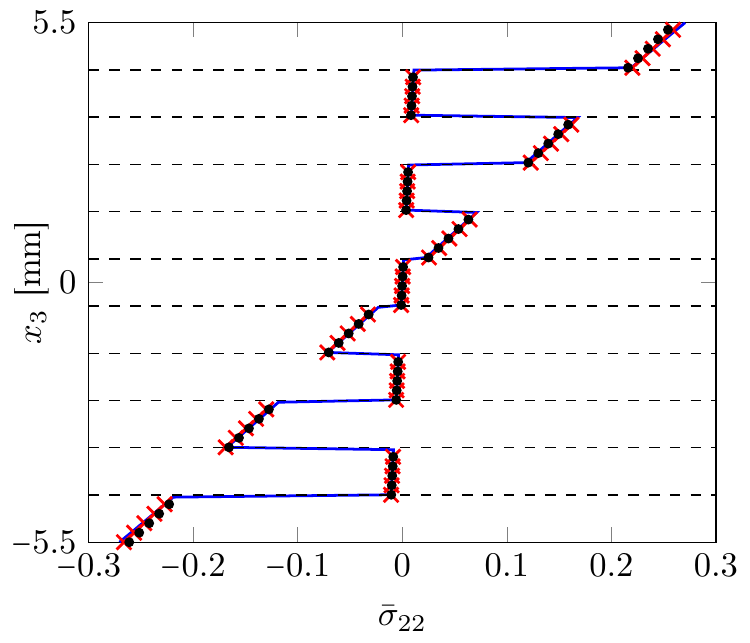}\fi}\hfill
	\subfigure[Normalized $\sigma_{22}$\label{subfig-4:in_plane_sig_1_1_l34}]{\ifrecompiletikz\tikzsetnextfilename{fig_03_b}\tikzexternalenable\input{images/fig_03_b}\tikzexternaldisable\else\includegraphics{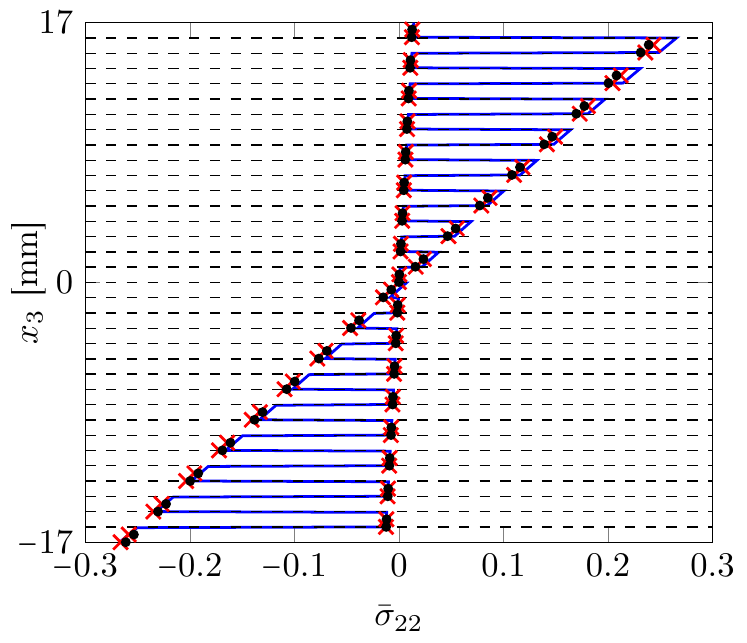}\fi}\\
	\subfigure[Normalized $\sigma_{12}$\label{subfig-5:in_plane_sig_1_1_l11}]{\ifrecompiletikz\tikzsetnextfilename{fig_02_c}\tikzexternalenable\input{images/fig_02_c}\tikzexternaldisable\else\includegraphics{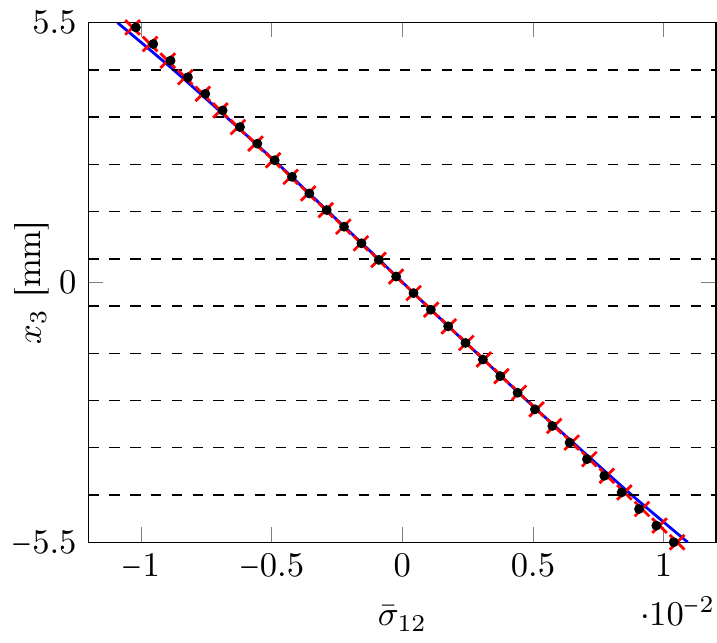}\fi}\hfill
	\subfigure[Normalized $\sigma_{12}$\label{subfig-6:in_plane_sig_1_1_l34}]{\ifrecompiletikz\tikzsetnextfilename{fig_03_c}\tikzexternalenable\input{images/fig_03_c}\tikzexternaldisable\else\includegraphics{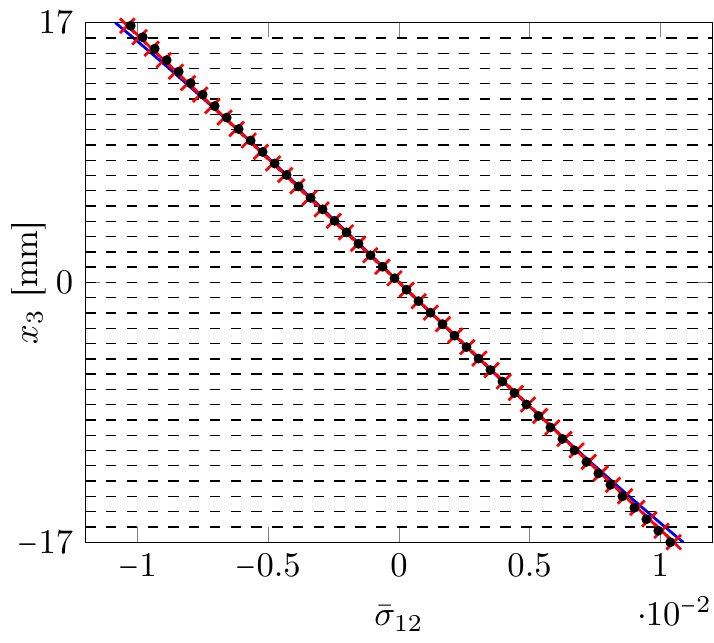}\fi}\
	\caption{Through-the-thickness in-plane stress solution for the Pagano problem~\cite{Pagano1970} evaluated at $x_1=x_2=L/4$. Plate case with 11 (left column) and 34 layers (right column), and length-to-thickness ratio $S=20$ (\blueline~Pagano's analytical solution versus numerical solutions obtained with $p=q=6$, and 7x7 control points:  \redcrosses~IGA-Galerkin, \blackcircles~IGA-Collocation).}
	\label{fig:example_in_plane_sig_1_1}
\end{figure}

As an example, in Figure \ref{fig:example_in_plane_sig_1_1} we present the in-plane solution profiles for a sampling point located at $x_1=x_2=L/4$, computed with both approaches described in Section \ref{subsec:collocation} and \ref{subsec:galerkin}, which prove to be accurate even for a rather small length-to-thickness plate ratio ($S=20$).

\begin{figure}[!htbp]
	\centering
	\subfigure[Normalized $\sigma_{13}$\label{subfig-1:out_of_plane_sig_1_1_l11}]{\ifrecompiletikz\tikzsetnextfilename{fig_02_d}\tikzexternalenable\input{images/fig_02_d}\tikzexternaldisable\else\includegraphics{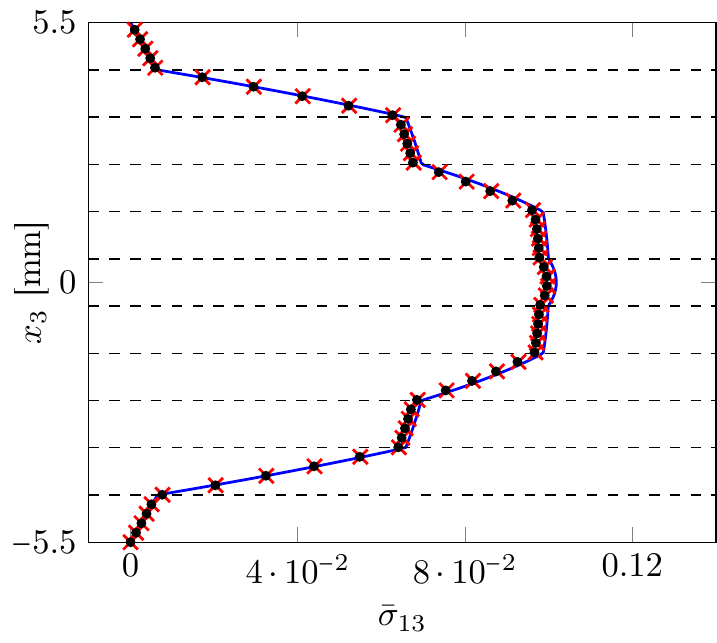}\fi}\hfill
	\subfigure[Normalized $\sigma_{13}$\label{subfig-2:out_of_plane_sig_1_1_l34}]{\ifrecompiletikz\tikzsetnextfilename{fig_03_d}\tikzexternalenable\input{images/fig_03_d}\tikzexternaldisable\else\includegraphics{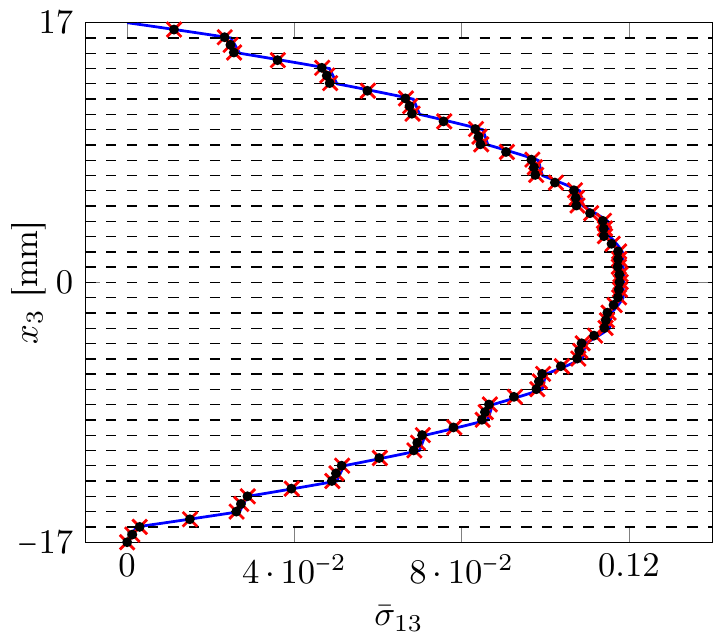}\fi}\\
	\subfigure[Normalized $\sigma_{23}$\label{subfig-2:out_of_plane_sig_1_1_lt11}]{\ifrecompiletikz\tikzsetnextfilename{fig_02_e}\tikzexternalenable\input{images/fig_02_e}\tikzexternaldisable\else\includegraphics{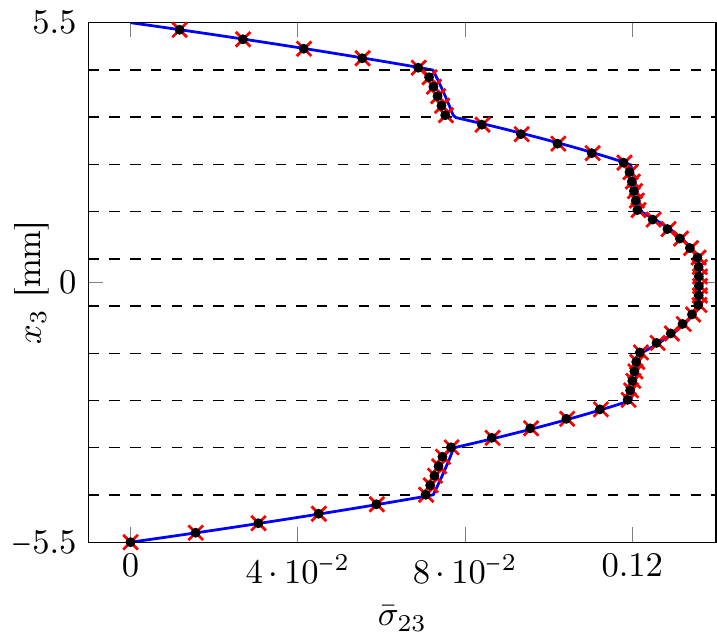}\fi}\hfill
	\subfigure[Normalized $\sigma_{23}$\label{subfig-4:out_of_plane_sig_1_1_l34}]{\ifrecompiletikz\tikzsetnextfilename{fig_03_e}\tikzexternalenable\input{images/fig_03_e}\tikzexternaldisable\else\includegraphics{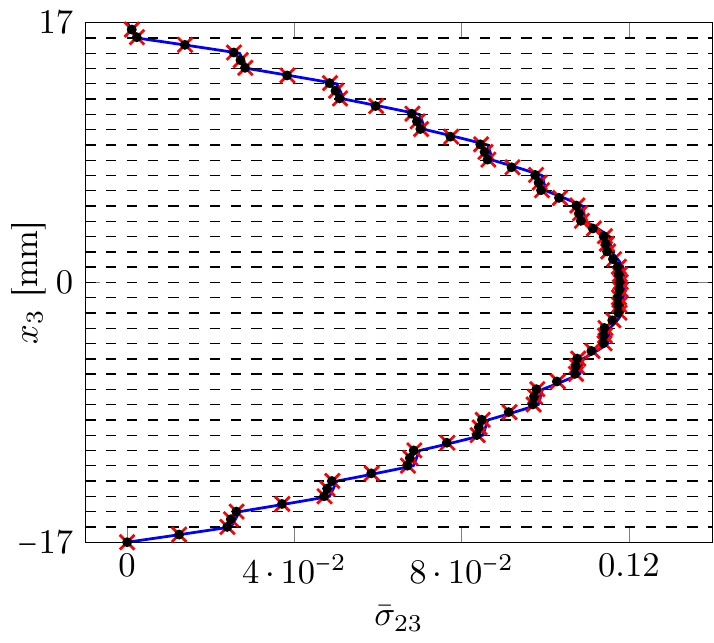}\fi}\\
	\subfigure[Normalized $\sigma_{33}$\label{subfig-5:out_of_plane_sig_1_1_l11}]{\ifrecompiletikz\tikzsetnextfilename{fig_02_f}\tikzexternalenable\input{images/fig_02_f}\tikzexternaldisable\else\includegraphics{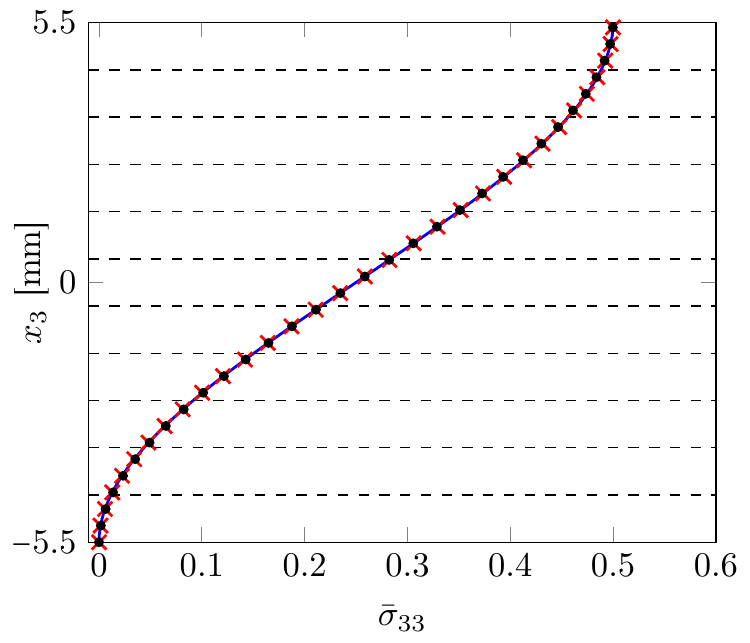}\fi}\hfill	
	\subfigure[Normalized $\sigma_{33}$\label{subfig-6:out_of_plane_sig_1_1_l34}]{\ifrecompiletikz\tikzsetnextfilename{fig_03_f}\tikzexternalenable\input{images/fig_03_f}\tikzexternaldisable\else\includegraphics{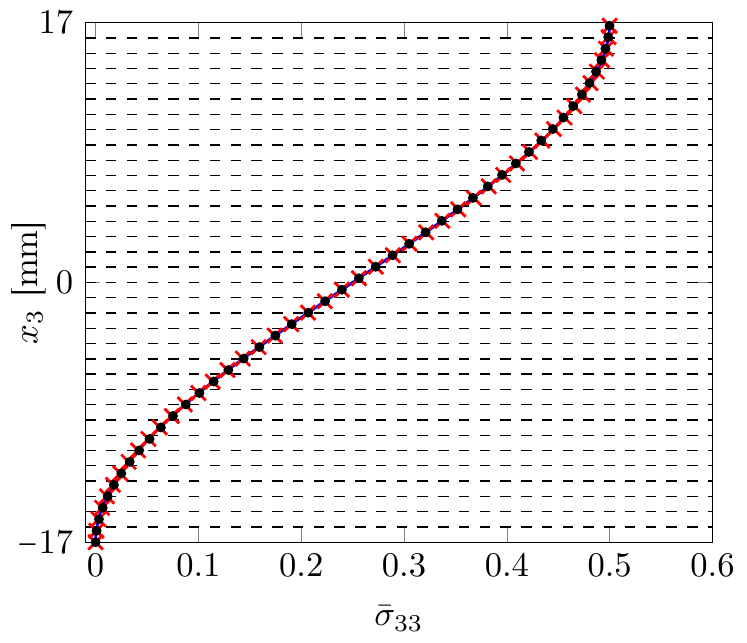}\fi}
	\caption{Through-the-thickness recovered out-of-plane stress solution for the Pagano problem~\cite{Pagano1970} evaluated at $x_1=x_2=L/4$. Plate case with 11 (left column) and 34 layers (right column), and length-to-thickness ratio $S=20$ (\blueline~Pagano's analytical solution versus post-processed numerical solutions obtained with degree of approximation $p=q=6$, and 7x7 control points:  \redcrosses~IGA-Galerkin, \blackcircles~IGA-Collocation).}
	\label{fig:example_out_of_plane_sig_1_1}
\end{figure}

In Figure \ref{fig:example_out_of_plane_sig_1_1} (for the same sampling point and plate geometrical description considered in Figure \ref{fig:example_in_plane_sig_1_1}) instead we readily reconstruct also an accurate out-of-plane stress state, applying the presented post-processing step based on equilibrium which can be regarded as inexpensive with respect to a full 3D analysis and is to be performed only at locations of interest.

We remark that using the CLPT to rigorously model non-symmetric cross-ply laminates, we would need to account for bending-stretching contributions.
However, for these type of laminates the bending-stretching coefficient matrix is not full and in addition the coupling effect decreases as the number of layers is increased \cite{Reddy2004,Whitney1996}.
Also, the presented numerical results are compared to Pagano's analytical solution, which is sufficiently general to describe the exact elastic response of rectangular, pinned edge laminates consisting of any number of orthotropic layers \cite{Pagano1970}.
Therefore, we can regard these coupling effects to be negligible and assume the proposed modeling approach to be an effective tool in understanding the behavior of the considered laminate class.
\begin{remark}\label{remark:1}
In \cite{Patton2019}, to tackle solid laminates, the presented equilibrium-based stress recovery procedure was combined, for collocation, with a homogenized through-the-thickness single-element approach, which is directly effective only for symmetric cross-ply distributions, as for non-symmetric ones the plate mid-plane is not balanced. Nevertheless, using the CLPT together with the proposed post-processing technique, we are able to accurately capture the behavior of non-symmetric cross-ply laminated plates also via collocation, despite neglecting bending-stretching coupling effects, and directly reconstruct the out-of-plane stresses from 2D displacement-based computations.
\end{remark}

\begin{figure}[!htbp]
\centering
\ifrecompiletikz\tikzsetnextfilename{fig_04}\tikzexternalenable\input{images/fig_04}\tikzexternaldisable\else\includegraphics{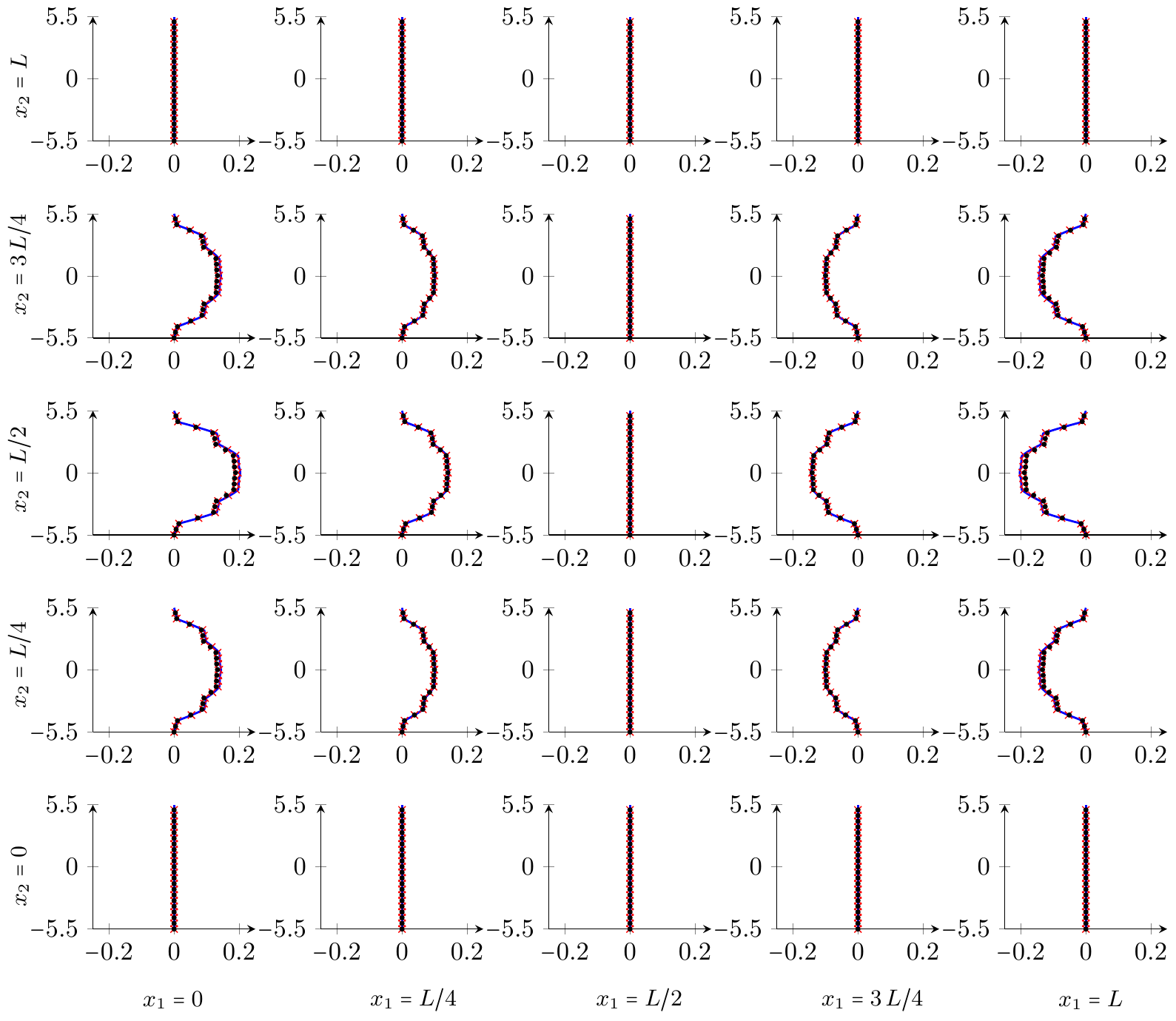}\fi
\caption{Through-the-thickness $\bar{\sigma}_{13}$ profiles for several in plane sampling points. $L$ represents the total length of the plate, that for this case is  $L=220\,\text{mm}$ (being $L=S\,t$ with $t=11\,\text{mm}$ and $S=20$), while the number of layers is 11 (\blueline~Pagano's analytical solution \cite{Pagano1970} versus recovered numerical solutions obtained with degree of approximation $p=q=6$, and 7x7 control points: \redcrosses~IGA-Galerkin, \blackcircles~IGA-Collocation).}
\label{fig:samplingS13_11lays}
\end{figure}

\begin{figure}[!htbp]
\centering
\ifrecompiletikz\tikzsetnextfilename{fig_05}\tikzexternalenable\input{images/fig_05}\tikzexternaldisable\else\includegraphics{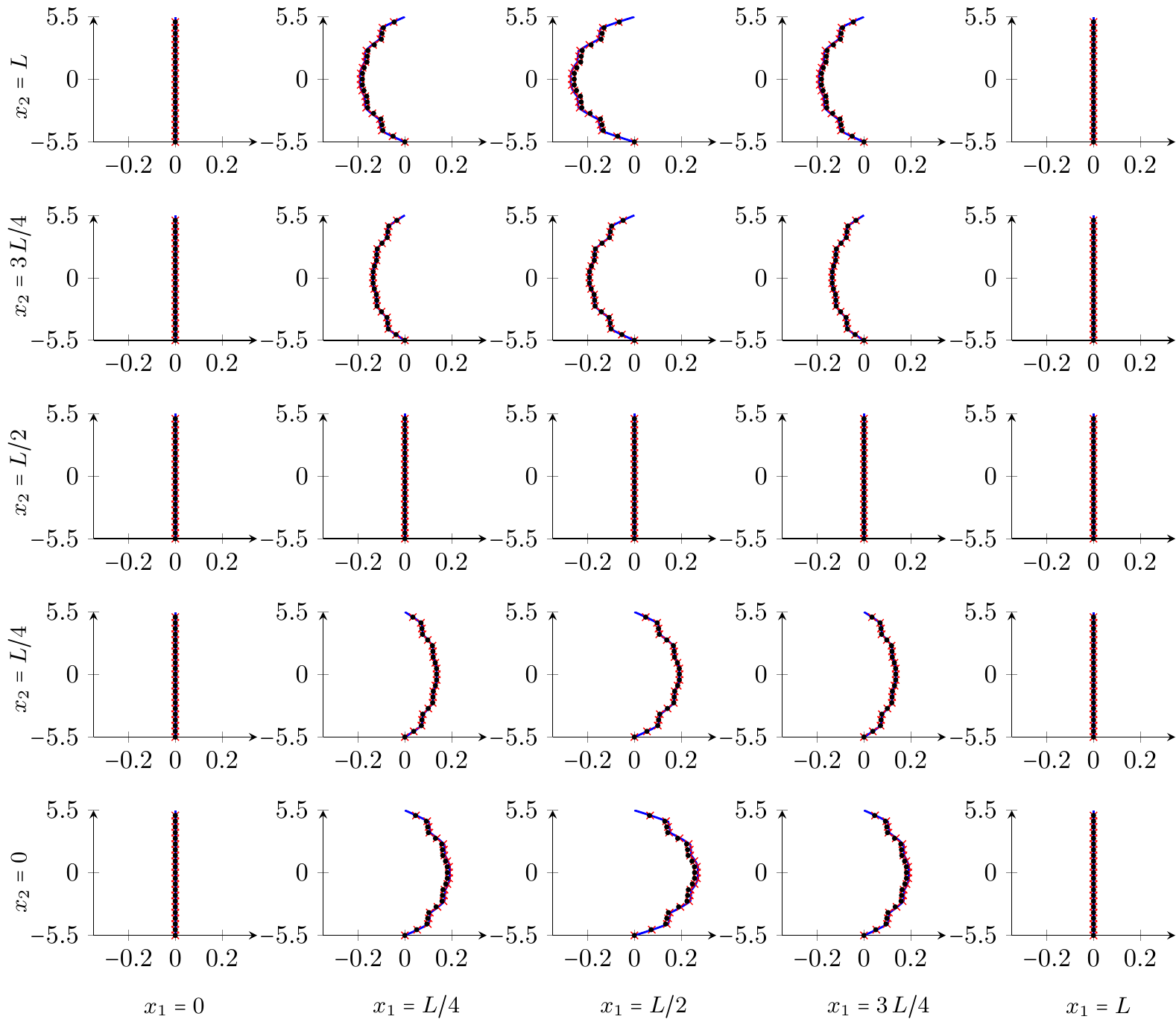}\fi
\caption{Through-the-thickness $\bar{\sigma}_{23}$ profiles for several in plane sampling points. $L$ represents the total length of the plate, that for this case is  $L=220\,\text{mm}$ (being $L=S\,t$ with $t=11\,\text{mm}$ and $S=20$), while the number of layers is 11 (\blueline~Pagano's analytical solution \cite{Pagano1970} versus recovered numerical solutions obtained with degree of approximation $p=q=6$, and 7x7 control points: \redcrosses~IGA-Galerkin, \blackcircles~IGA-Collocation).}
\label{fig:samplingS23_11lays}
\end{figure}

\begin{figure}[!htbp]
\centering
\ifrecompiletikz\tikzsetnextfilename{fig_06}\tikzexternalenable\input{images/fig_06}\tikzexternaldisable\else\includegraphics{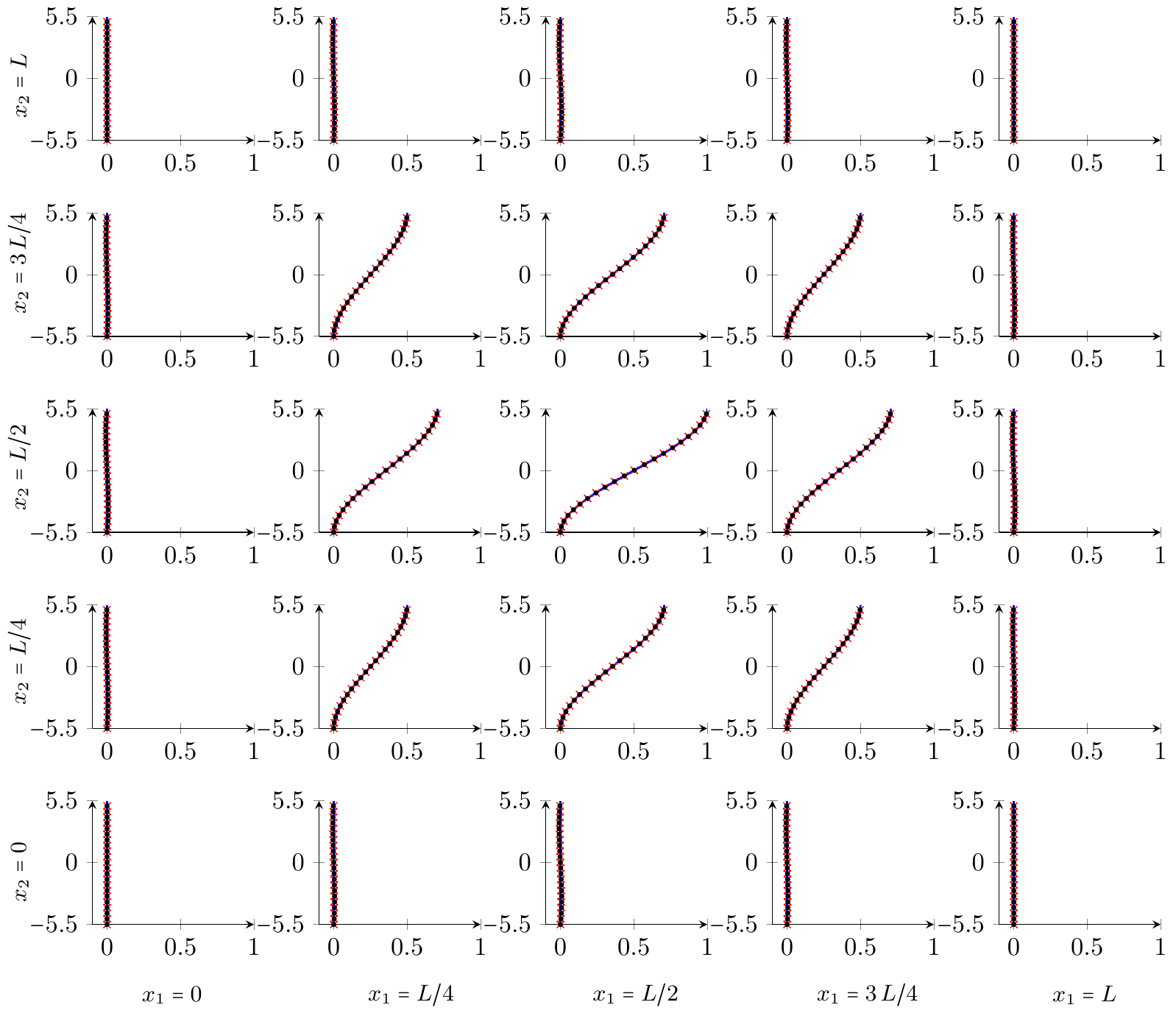}\fi
\caption{Through-the-thickness $\bar{\sigma}_{33}$ profiles for several in plane sampling points. $L$ represents the total length of the plate, that for this case is  $L=220\,\text{mm}$ (being $L=S\,t$ with $t=11\,\text{mm}$ and $S=20$), while the number of layers is 11 (\blueline~Pagano's analytical solution \cite{Pagano1970} versus recovered numerical solutions obtained with degree of approximation $p=q=6$, and 7x7 control points: \redcrosses~IGA-Galerkin, \blackcircles~IGA-Collocation).}
\label{fig:samplingS33_11lays}
\end{figure}

\begin{figure}[!htbp]
	\centering
	\ifrecompiletikz\tikzsetnextfilename{fig_07}\tikzexternalenable\input{images/fig_07}\tikzexternaldisable\else\includegraphics{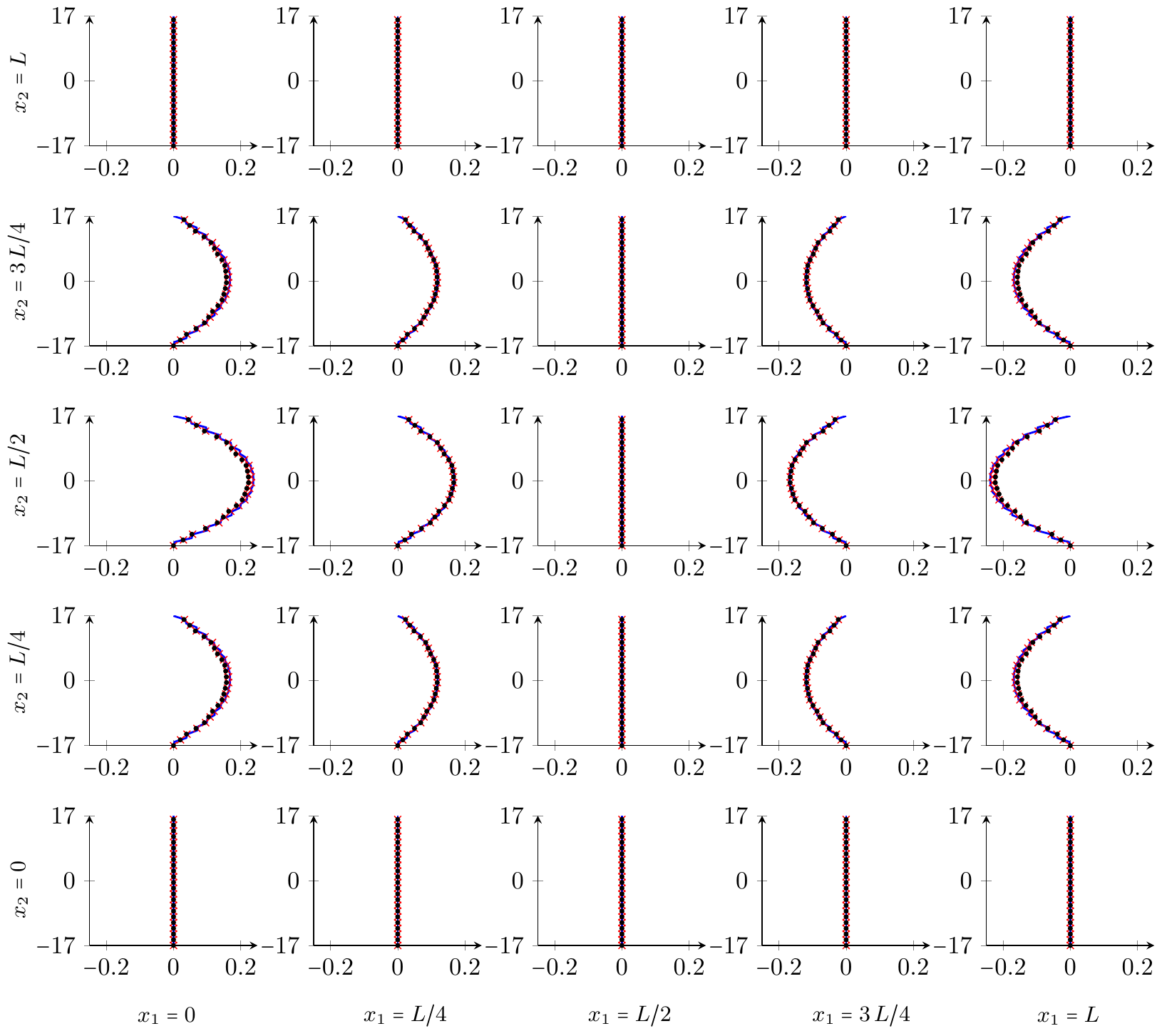}\fi
	\caption{Through-the-thickness $\bar{\sigma}_{13}$ profiles for several in plane sampling points. $L$ represents the total length of the plate, that for this case is  $L=1020\,\text{mm}$ (being $L=S\,t$ with $t=34\,\text{mm}$ and $S=20$), while the number of layers is 34 (\blueline~Pagano's analytical solution \cite{Pagano1970} versus recovered numerical solutions obtained with degree of approximation $p=q=6$, and 7x7 control points: \redcrosses~IGA-Galerkin, \blackcircles~IGA-Collocation).}
	\label{fig:samplingS13_34lays}
\end{figure}

\begin{figure}[!htbp]
	\centering
	\ifrecompiletikz\tikzsetnextfilename{fig_08}\tikzexternalenable\input{images/fig_08}\tikzexternaldisable\else\includegraphics{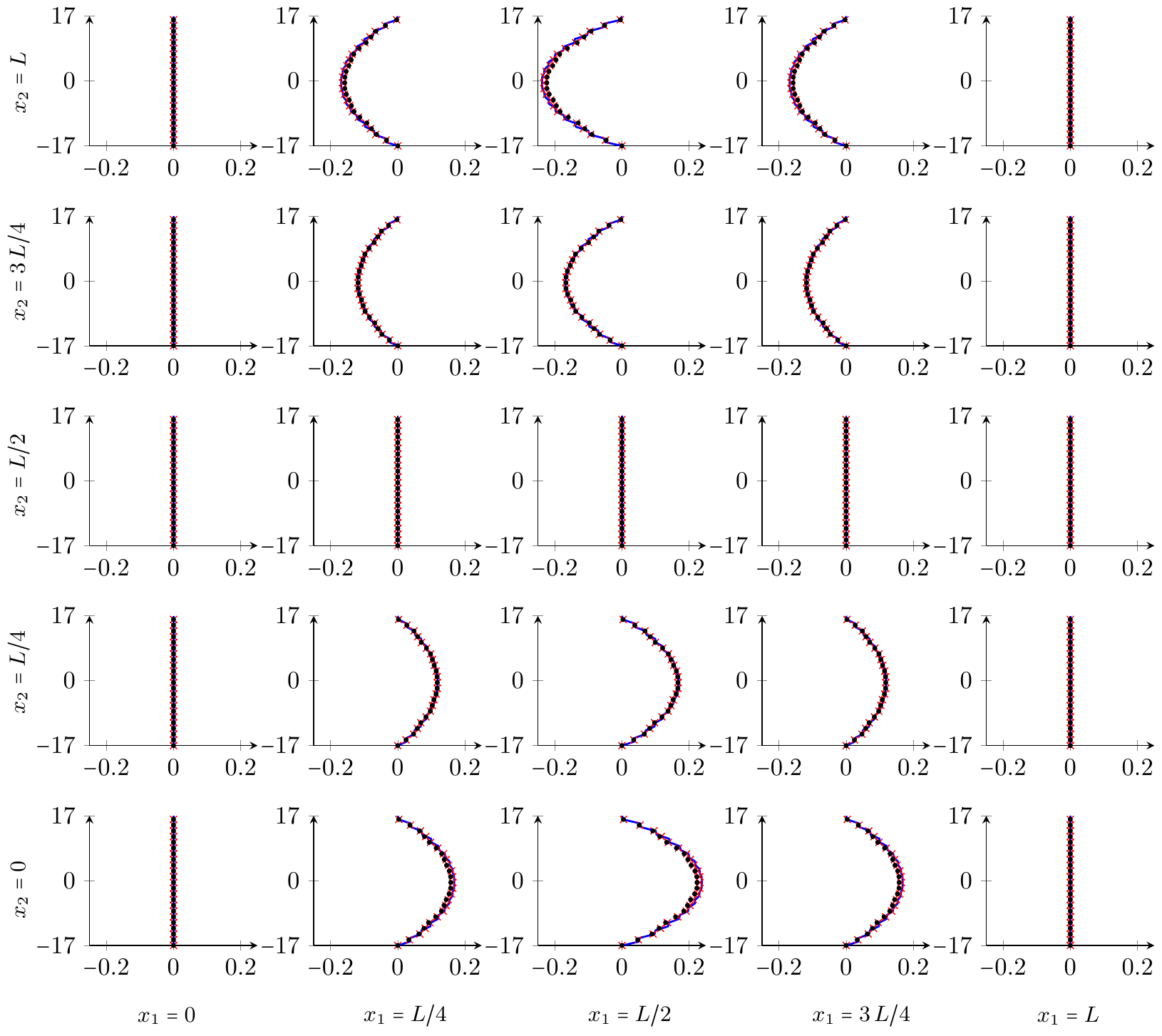}\fi
	\caption{Through-the-thickness $\bar{\sigma}_{23}$ profiles for several in plane sampling points. $L$ represents the total length of the plate, that for this case is  $L=1020\,\text{mm}$ (being $L=S\,t$ with $t=34\,\text{mm}$ and $S=20$), while the number of layers is 34 (\blueline~Pagano's analytical solution \cite{Pagano1970} versus recovered numerical solutions obtained with degree of approximation $p=q=6$, and 7x7 control points: \redcrosses~IGA-Galerkin, \blackcircles~IGA-Collocation).}
	\label{fig:samplingS23_34lays}
\end{figure}

\begin{figure}[!htbp]
	\centering
	\ifrecompiletikz\tikzsetnextfilename{fig_09}\tikzexternalenable\input{images/fig_09}\tikzexternaldisable\else\includegraphics{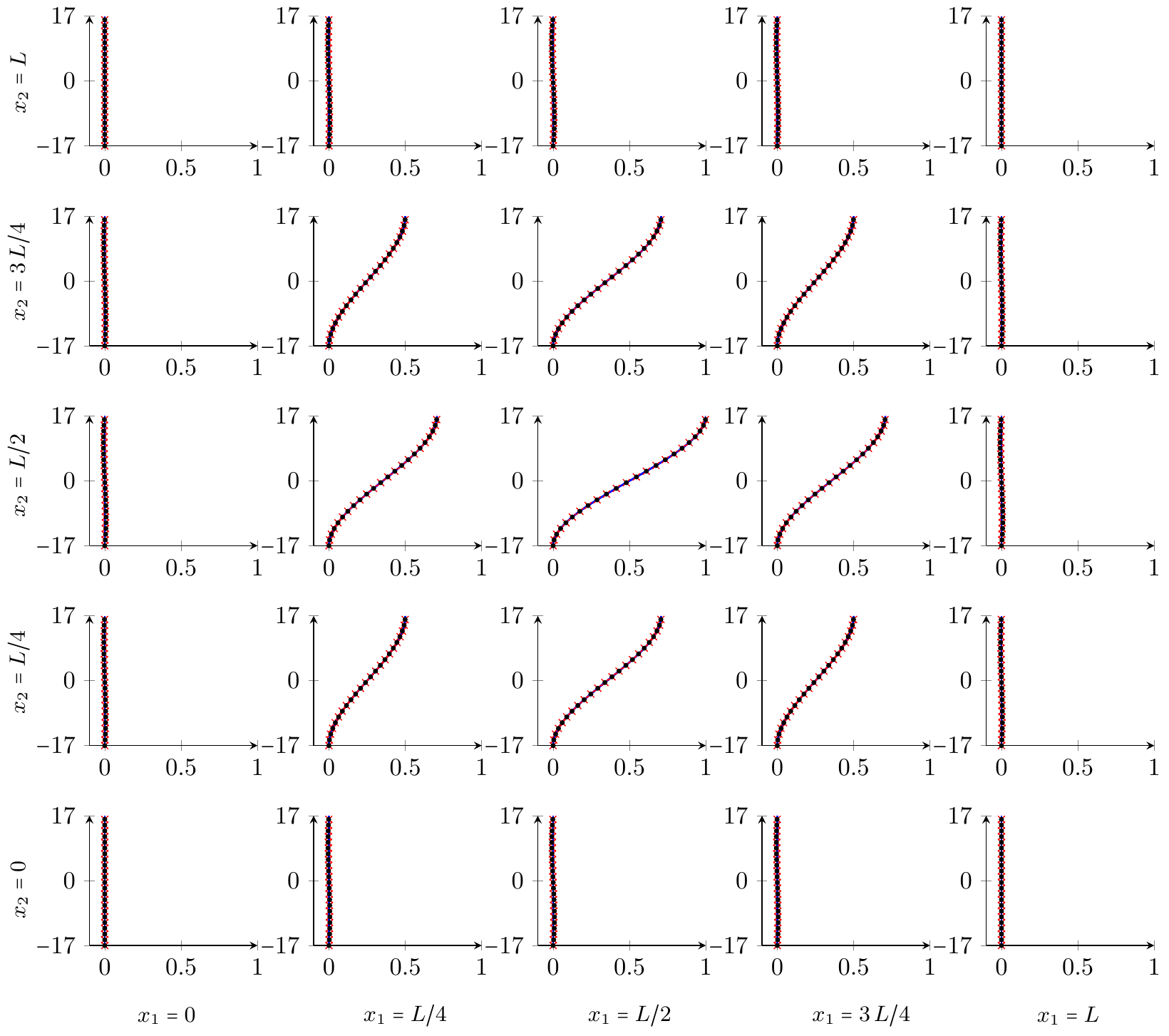}\fi
	\caption{Through-the-thickness $\bar{\sigma}_{33}$ profiles for several in plane sampling points. $L$ represents the total length of the plate, that for this case is  $L=1020\,\text{mm}$ (being $L=S\,t$ with $t=34\,\text{mm}$ and $S=20$), while the number of layers is 34 (\blueline~Pagano's analytical solution \cite{Pagano1970} versus recovered numerical solutions obtained with degree of approximation $p=q=6$, and 7x7 control points: \redcrosses~IGA-Galerkin, \blackcircles~IGA-Collocation).}
	\label{fig:samplingS33_34lays}
\end{figure}

\newpage
In Figures \ref{fig:samplingS13_11lays}-\ref{fig:samplingS33_34lays} the out-of-plane stress state profile is recovered
sampling the composite plate every quarter of length in both in-plane directions, to show the effect of post-processing at different locations of the plate for both a symmetric and a non-symmetric ply distribution of 11 and 34 layers, respectively (see Remark \ref{remark:1}).
Across all sampled points, the proposed approach accurately captures the 3D stresses in every single
layer when compared to Pagano's solution. Also, the model remains accurate at the boundaries, where solution inaccuracy is typically expected \cite{Mittelstedt2007}, and satisfies the traction-free conditions for transverse shear stresses at the top and bottom surfaces of the laminate.
\newpage

\subsection{Parametric study on length-to-thickness ratio}\label{subsec:param_study}
In order to further investigate the proposed approach, examples which consider a varying length-to-thickness ratio (i.e., $S=20,30,40,\;\text{and}\;50$) are performed respectively
for 11 and 34 layers, examining an increasing number of degrees of freedom.

\begin{figure}[!htbp]
	\centering
	\hspace{-7pt}
	\subfigure[11 layers \label{subfig-4:errors13_l11_6_varSGal}]{\ifrecompiletikz\tikzsetnextfilename{fig_10_d}\tikzexternalenable\input{images/fig_10_d}\tikzexternaldisable\else\includegraphics{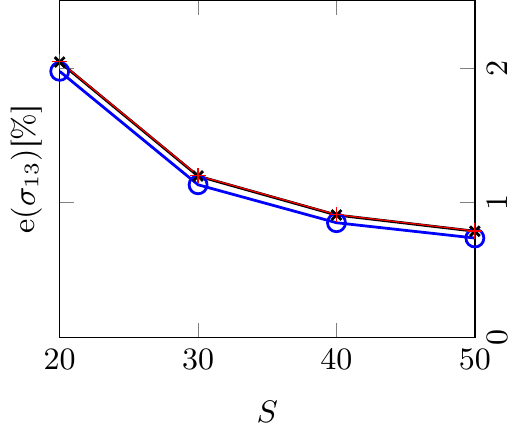}\fi}\hspace{-7pt}
	\subfigure[11 layers \label{subfig-5:errors23_l11_6_varSGal}]{\ifrecompiletikz\tikzsetnextfilename{fig_10_e}\tikzexternalenable\input{images/fig_10_e}\tikzexternaldisable\else\includegraphics{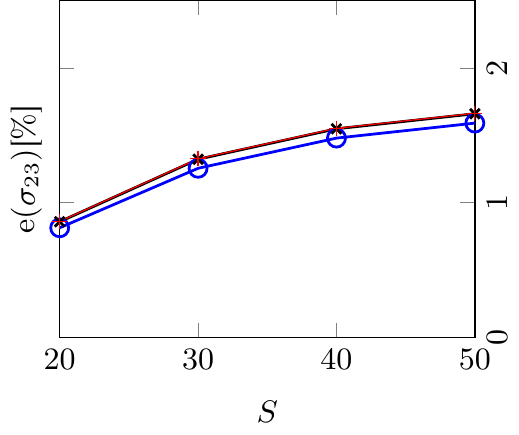}\fi}\hspace{-7pt}
	\subfigure[11 layers \label{subfig-6:errors33_l11_6_varSGal}]{\ifrecompiletikz\tikzsetnextfilename{fig_10_f}\tikzexternalenable\input{images/fig_10_f}\tikzexternaldisable\else\includegraphics{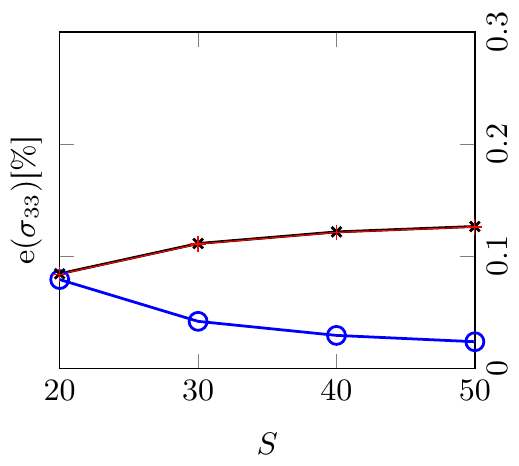}\fi}\\
	\hspace{-7pt}
	\subfigure[34 layers \label{subfig-7:errors13_l34_6_varSGal}]{\ifrecompiletikz\tikzsetnextfilename{fig_11_d}\tikzexternalenable\input{images/fig_11_d}\tikzexternaldisable\else\includegraphics{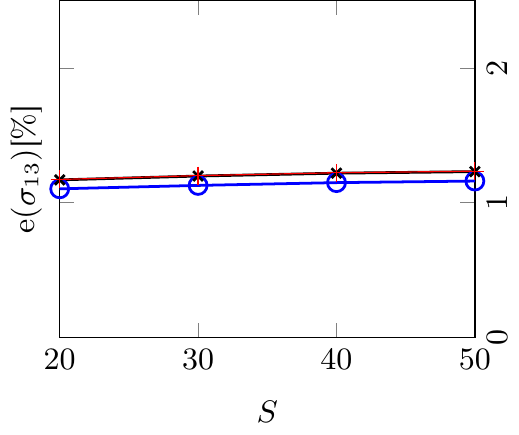}\fi}\hspace{-7pt}
	\subfigure[34 layers \label{subfig-8:errors23_l34_6_varSGal}]{\ifrecompiletikz\tikzsetnextfilename{fig_11_e}\tikzexternalenable\input{images/fig_11_e}\tikzexternaldisable\else\includegraphics{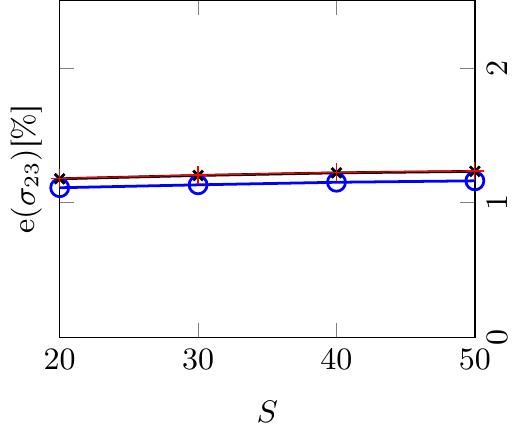}\fi}\hspace{-7pt}
	\subfigure[34 layers \label{subfig-9:errors33_l34_6_varSGal}]{\ifrecompiletikz\tikzsetnextfilename{fig_11_f}\tikzexternalenable\input{images/fig_11_f}\tikzexternaldisable\else\includegraphics{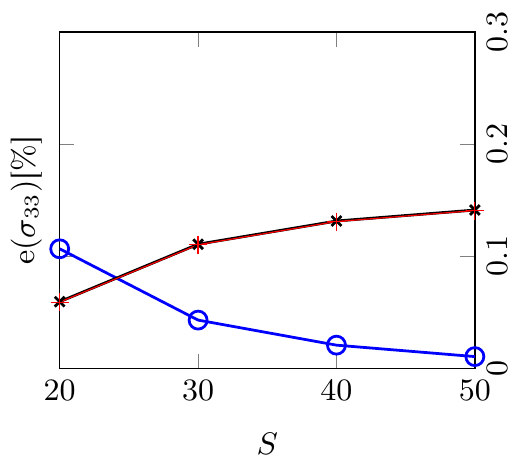}\fi}
	\caption{$L^2$ relative percentage error evaluation at $x_1=x_2=L/4$ for IGA-Galerkin using an in-plane degree of approximation equal to 6. Different length-to-thickness ratios $S$ are investigated for a number of layers equal to 11 and 34 (Number of control points per in-plane direction: 7 \bluesolidcircle, 14 \blacksolidx, 21 \redsolidcross).}
	\label{fig:conv6_varSGal}
\end{figure}

\begin{figure}[!htbp]
	\centering
	\hspace{-7pt}
	\subfigure[11 layers \label{subfig-4:errors13_l11_6_varSColl}]{\ifrecompiletikz\tikzsetnextfilename{fig_12_d}\tikzexternalenable\input{images/fig_12_d}\tikzexternaldisable\else\includegraphics{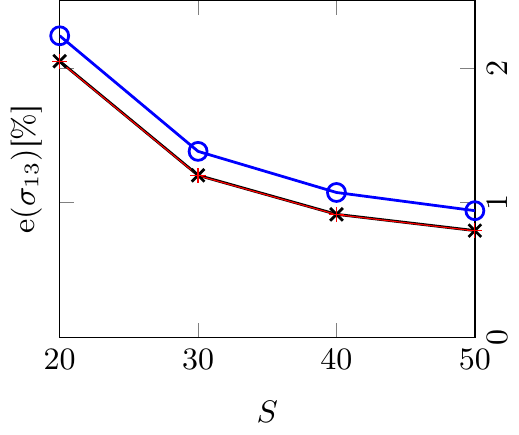}\fi}\hspace{-7pt}
	\subfigure[11 layers \label{subfig-5:errors23_l11_6_varSColl}]{\ifrecompiletikz\tikzsetnextfilename{fig_12_e}\tikzexternalenable\input{images/fig_12_e}\tikzexternaldisable\else\includegraphics{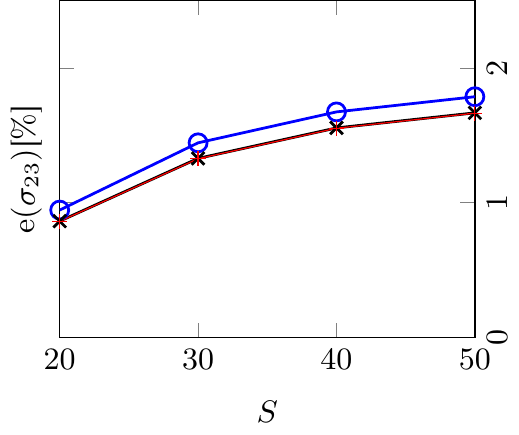}\fi}\hspace{-7pt}
	\subfigure[11 layers \label{subfig-6:errors33_l11_6_varSColl}]{\ifrecompiletikz\tikzsetnextfilename{fig_12_f}\tikzexternalenable\input{images/fig_12_f}\tikzexternaldisable\else\includegraphics{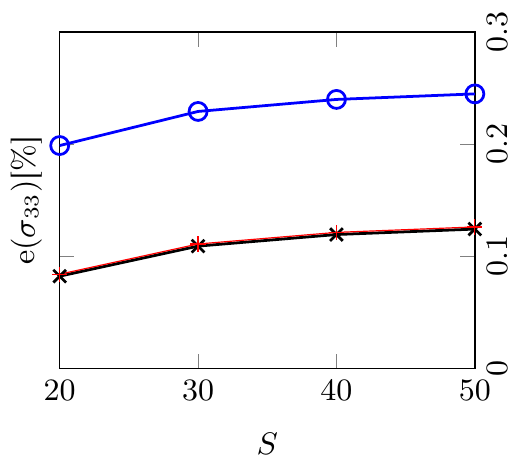}\fi}\\
	\hspace{-7pt}
	\subfigure[34 layers \label{subfig-7:errors13_l34_6_varSColl}]{\ifrecompiletikz\tikzsetnextfilename{fig_13_d}\tikzexternalenable\input{images/fig_13_d}\tikzexternaldisable\else\includegraphics{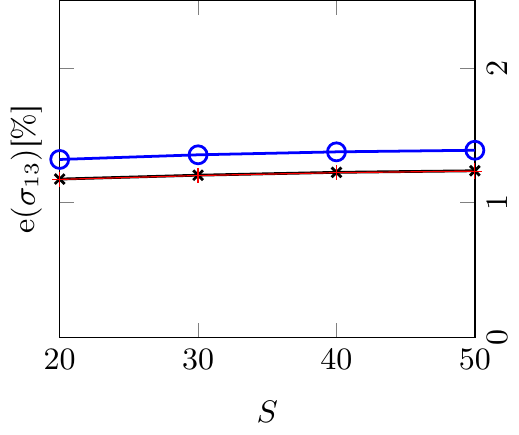}\fi}\hspace{-7pt}
	\subfigure[34 layers \label{subfig-8:errors23_l34_6_varSColl}]{\ifrecompiletikz\tikzsetnextfilename{fig_13_e}\tikzexternalenable\input{images/fig_13_e}\tikzexternaldisable\else\includegraphics{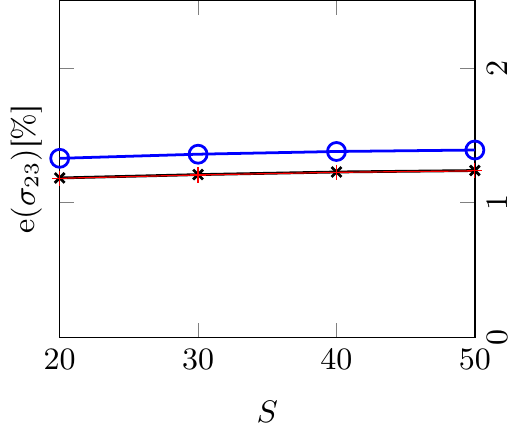}\fi}\hspace{-7pt}
	\subfigure[34 layers \label{subfig-9:errors33_l34_6_varSColl}]{\ifrecompiletikz\tikzsetnextfilename{fig_13_f}\tikzexternalenable\input{images/fig_13_f}\tikzexternaldisable\else\includegraphics{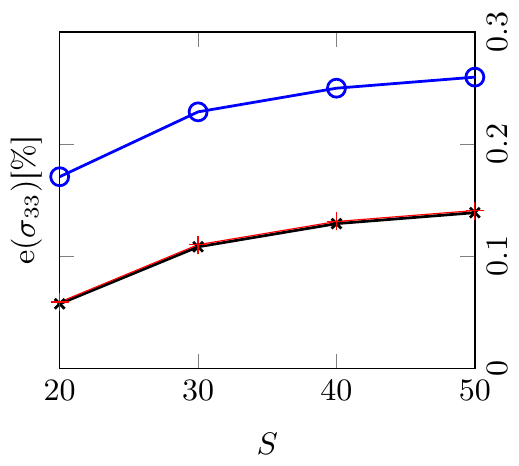}\fi}
	\caption{$L^2$ relative percentage error evaluation at $x_1=x_2=L/4$ for IGA-Collocation using an in-plane degree of approximation equal to 6. Different length-to-thickness ratios $S$ are investigated for a number of layers equal to 11 and 34 (Number of control points per in-plane direction: 7 \bluesolidcircle, 14 \blacksolidx, 21 \redsolidcross).}
	\label{fig:conv6_varSColl}
\end{figure}

In Figure \ref{fig:conv6_varSGal} and \ref{fig:conv6_varSColl} we assess the performance of both the isogeometric Galerkin and the collocation approach coupled with the presented post-processing technique at $x_1=x_2=L/4$, adopting the following $L^2$ error definition
\begin{alignat}{2}\label{eq:error}
&\text{e}(\sigma_{i3})=\sqrt{\cfrac{\int_{x_3}(\sigma_{i3}^\text{analytic}(\bar{x}_1,\bar{x}_2,x_3)-\sigma_{i3}^\text{recovered}(\bar{x}_1,\bar{x}_2,x_3))^2}{\int_{x_3}(\sigma_{i3}^\text{analytic}(\bar{x}_1,\bar{x}_2,x_3))^2}} \quad&&\quad i=1,2,3\,.
\end{alignat}
The post-processing approach seems to be particularly
suitable for tackling plates characterized by a significant number of layers.
Moreover, we observe that the modeling error, given by the a-posteriori step, dominates over the approximation one; thus, further refinement operations do not seem to provide a significant benefit for the considered tests. We want to highlight, however, that errors are typically in the 1.5\% range or lower in this case.

We would like to remark that further tests have been carried out for a lower degree displacement field approximation (i.e., $4 \le p=q<6$), which led to a less accurate out-of-plane stress reconstruction in particular for collocation. In our experience, adopting a degree of approximation equal to 6 seems to be a reasonable choice to correctly reproduce the complete 3D stress state for both considered methods. In such a case using only one element to approximate the plate mid-plane, corresponding to 49 d.o.f's, is sufficient to provide good results in the considered example, which is characterized by a simple geometry.

\newpage
\subsection{Assessment of the approach at the plate boundary}\label{sec:boundary}
Interlaminar stresses in laminates subjected to transverse loadings may become important near the structure edges. For example in the case of cross-ply laminates out-of-plane
stresses usually face weaker material strength properties according to the stacking sequence, leading in the proximity of material discontinuities to
stress concentrations, which may result in
premature failure of the structure due to delamination
fracture \cite{Mittelstedt2007}.

Thus, we further test the proposed post-processing technique, studying both symmetric and non-symmetric cross-ply plate cases (namely, stacking sequences made of 11 and 34 layers) investigating the composite behavior especially at the boundary. To this extent, we consider an increasing length-to-thickness ratio (i.e., $S= 20,\;30,\;40,\;\text{and}\;50$) for fixed degrees of approximation $p = q = 6$ using 7x7 control points, and we report in Tables  \ref{tab:error34layers_pnt1}-\ref{tab:error34layers_pnt5} the out-of-plane stress pointwise relative difference defined as
\begin{alignat}{2}\label{eq:diffrel}
&\Delta(\sigma_{i3})=\cfrac{\vert\sigma_{i3}^\text{analytic}(\bar{x}_1,\bar{x}_2,\bar{x}_3)-\sigma_{i3}^\text{recovered}(\bar{x}_1,\bar{x}_2,\bar{x}_3)\vert}{\vert\sigma_{i3}^\text{analytic}(\bar{x}_1,\bar{x}_2,\bar{x}_3)\vert}\quad&&\quad i=1,2,3\,.
\end{alignat}
Note that to avoid divisions by zero, when the corresponding analytical solution is zero we compute difference values marked with the * symbol instead, computed as  
\begin{alignat}{2}\label{eq:diffabs}
&\Delta(\sigma_{i3})=\vert\sigma_{i3}^\text{analytic}(\bar{x}_1,\bar{x}_2,\bar{x}_3)-\sigma_{i3}^\text{recovered}(\bar{x}_1,\bar{x}_2,\bar{x}_3)\vert\quad&&\quad i=1,2,3\,.
\end{alignat}
In \eqref{eq:diffrel} and \eqref{eq:diffabs} the $\bar{.}$ symbol means that a fixed coordinate in the plate domain is considered.

For the selected sampling points, a single in-plane element comprising 7x7 degrees of freedom is able to provide for a plate made of 11 layers, maximum differences of 4\% or lower (3\% or lower for a 34-layered plate) on the boundary and of 2.5\% or lower (less than 1\% for a 34-layered case) inside the domain for the considered isogeometric Galerkin method coupled with the proposed post-processing technique. 
Under the same modeling conditions, collocation combined with the equilibrium-based strategy allows to obtain maximum differences of 8\% or lower (6.5\% or lower for a 34-layered plate) on the border and of 3\% or lower (1.5\% or lower for a 34-layered case) inside the plate. 
Finally, relative differences with reference to normal out-of-plane $\sigma_{33}$ are, on average, one order magnitude less than those that correspond to the shear stress counterparts.

\begin{table}[!htbp]\caption{Simply supported composite plate under a sinusoidal load with 11 layers. Out-of-plane stress state difference with respect to Pagano's solution~\cite{Pagano1970}. We compare, at $\boldsymbol{x}=(0,L/2,0)$, post-processed isogeometric collocation approach (IGA-C) and post-processed isogeometric Galerkin method (IGA-G) for a degree of approximation $p=q=6$ and 7x7 control points.}
	\vspace{0.5cm}
	\centering	
	\begin{adjustbox}{max width=\textwidth}
		\begin{tabular}{? c ? c ? c | c | c ? c | c | c ?}
			\thickhline
			\multirow{2}{*}{$\boldsymbol{S}$}&\multirow{2}{*}{\textbf{Method}}&$\sigma_{13}(0,L/2,0)$
			&$\sigma_{23}(0,L/2,0)$
			&$\sigma_{33}(0,L/2,0)$& $\Delta(\sigma_{13})$&$\Delta(\sigma_{23})$&$\Delta(\sigma_{33})$ \Tstrut\Bstrut\\
			&&[-]&[-]&[-]&[\%]&[\%]&[\%]\Tstrut\Bstrut\\\thickhline
			\multirow{3}{*}{20}&Analytical& 4.0728&0.0000&0.0000&-&-&-
			\Tstrut\Bstrut\\\cline{2-8}
			&post-processed IGA-G &3.9290&0.0000&0.0001&3.5295&0.0000*&0.0055*
			\Tstrut\Bstrut\\\cline{2-8}
			&post-processed IGA-C &3.7848&0.0000&0.0001&7.0706&0.0000*& 0.0100*\Tstrut\Bstrut\\\thickhline
			\multirow{3}{*}{30}&Analytical& 6.0598&0.0000&0.0000&-&-&-
			\Tstrut\Bstrut\\\cline{2-8}
			&post-processed IGA-G &5.8935&0.0000&0.0001&2.7445&0.0000*&0.0055*
			\Tstrut\Bstrut\\\cline{2-8}
			&post-processed IGA-C &5.6772&0.0000&0.0001&6.3145&0.0000*&0.0100*\Tstrut\Bstrut\\\thickhline
			\multirow{3}{*}{40}&Analytical& 8.0545&0.0000&0.0000&-&-&-
			\Tstrut\Bstrut\\\cline{2-8}
			&post-processed IGA-G &7.8580&0.0000&0.0001&2.4395&0.0000*&0.0055*
			\Tstrut\Bstrut\\\cline{2-8}
			&post-processed IGA-C &7.5696&0.0000&0.0001&6.0206&0.0000*& 0.0100*\Tstrut\Bstrut\\\thickhline
			\multirow{3}{*}{50}&Analytical& 10.0530&0.0000&0.0000&-&-&-
			\Tstrut\Bstrut\\\cline{2-8}
			&post-processed IGA-G &9.8225&0.0000&-0.0001&2.2923&0.0000*&0.0055*
			\Tstrut\Bstrut\\\cline{2-8}
			&post-processed IGA-C &9.4620&0.0000&0.0001&5.8789&0.0000*& 0.0100*\Tstrut\Bstrut\\
			\thickhline
		\end{tabular}
	\end{adjustbox}\label{tab:error11layers_pnt1}
\end{table}
\begin{table}[!htbp]\caption{Simply supported composite plate under a sinusoidal load with 11 layers. Out-of-plane stress state difference with respect to Pagano's solution~\cite{Pagano1970}. We compare, at $\boldsymbol{x}=(0,L/2,h/4)$, post-processed isogeometric collocation approach (IGA-C) and post-processed isogeometric Galerkin method (IGA-G) for a degree of approximation $p=q=6$ and 7x7 control points.}
	\vspace{0.5cm}
	\centering	
	\begin{adjustbox}{max width=\textwidth}
		\begin{tabular}{? c ? c ? c | c | c ? c | c | c ?}
			\thickhline
			\multirow{2}{*}{$\boldsymbol{S}$}&\multirow{2}{*}{\textbf{Method}}&$\sigma_{13}(0,L/2,h/4)$
			&$\sigma_{23}(0,L/2,h/4)$
			&$\sigma_{33}(0,L/2,h/4)$& $\Delta(\sigma_{13})$&$\Delta(\sigma_{23})$&$\Delta(\sigma_{33})$ \Tstrut\Bstrut\\
			&&[-]&[-]&[-]&[\%]&[\%]&[\%]\Tstrut\Bstrut\\\thickhline
			\multirow{3}{*}{20}&Analytical& 2.7527&0.0000&0.0000&-&-&-
			\Tstrut\Bstrut\\\cline{2-8}
			&post-processed IGA-G &2.6394&0.0000&-0.0045&4.1167&0.0000*&0.4464*
			\Tstrut\Bstrut\\\cline{2-8}
			&post-processed IGA-C &2.5433&0.0000&-0.0082&7.6056&0.0000*& 0.8161*\Tstrut\Bstrut\\\thickhline
			\multirow{3}{*}{30}&Analytical& 4.0817&0.0000&0.0000&-&-&-
			\Tstrut\Bstrut\\\cline{2-8}
			&post-processed IGA-G &3.9590&0.0000&-0.0045&3.0059&0.0000*&0.4464*
			\Tstrut\Bstrut\\\cline{2-8}
			&post-processed IGA-C &3.8150&0.0000&-0.0082&6.5352&0.0000*& 0.8161*\Tstrut\Bstrut\\\thickhline
			\multirow{3}{*}{40}&Analytical& 5.4188&0.0000&0.0000&-&-&-
			\Tstrut\Bstrut\\\cline{2-8}
			&post-processed IGA-G &5.2787&0.0000&-0.0045&2.5849&0.0000*&0.4464*
			\Tstrut\Bstrut\\\cline{2-8}
			&post-processed IGA-C &5.0866&0.0000&-0.0082&6.1295&0.0000*& 0.8161*\Tstrut\Bstrut\\\thickhline
			\multirow{3}{*}{50}&Analytical& 6.7595&0.0000&0.0000&-&-&-
			\Tstrut\Bstrut\\\cline{2-8}
			&post-processed IGA-G &6.5984&0.0000&-0.0045&2.3838&0.0000*&0.4464*
			\Tstrut\Bstrut\\\cline{2-8}
			&post-processed IGA-C &6.3583&0.0000&-0.0082&5.9357&0.0000*& 0.8161*\Tstrut\Bstrut\\
			\thickhline
		\end{tabular}
	\end{adjustbox}\label{tab:error11layers_pnt2}
\end{table}
\begin{table}[!htbp]\caption{Simply supported composite plate under a sinusoidal load with 11 layers. Out-of-plane stress state difference with respect to Pagano's solution~\cite{Pagano1970}. We compare, at $\boldsymbol{x}=(L/4,L/4,0)$, post-processed isogeometric collocation approach (IGA-C) and post-processed isogeometric Galerkin method (IGA-G) for a degree of approximation $p=q=6$ and 7x7 control points.}
	\vspace{0.5cm}
	\centering	
	\begin{adjustbox}{max width=\textwidth}
		\begin{tabular}{? c ? c ? c | c | c ? c | c | c ?}
			\thickhline
			\multirow{2}{*}{$\boldsymbol{S}$}&\multirow{2}{*}{\textbf{Method}}&$\sigma_{13}$ $(L/4,L/4,0)$
			&$\sigma_{23}$ $(L/4,L/4,0)$
			&$\sigma_{33}$ $(L/4,L/4,0)$& $\Delta(\sigma_{13})$&$\Delta(\sigma_{23})$&$\Delta(\sigma_{33})$ \Tstrut\Bstrut\\
			&&[-]&[-]&[-]&[\%]&[\%]&[\%]\Tstrut\Bstrut\\\thickhline
			\multirow{3}{*}{20}&Analytical&2.0364&2.7220&0.2483&-&-&-
			\Tstrut\Bstrut\\\cline{2-8}
			&post-processed IGA-G &1.9974&2.7240&0.2483&1.9166&0.0751&0.0026
			\Tstrut\Bstrut\\\cline{2-8}
			&post-processed IGA-C &1.9919&2.7183&0.2483&2.1852&0.1340&0.0080\Tstrut\Bstrut\\\thickhline
			\multirow{3}{*}{30}&Analytical& 3.0299&4.1212&0.2483&-&-&-
			\Tstrut\Bstrut\\\cline{2-8}
			&post-processed IGA-G &2.9960&4.0860&0.2483&1.1185&0.8526&0.0000
			\Tstrut\Bstrut\\\cline{2-8}
			&post-processed IGA-C &2.9878&4.0775&0.2483&1.3893&1.0598&0.0054\Tstrut\Bstrut\\\thickhline
			\multirow{3}{*}{40}&Analytical& 4.0273&5.5138&0.2483&-&-&-
			\Tstrut\Bstrut\\\cline{2-8}
			&post-processed IGA-G &3.9947&5.4481&0.2483&0.8083&1.1926&0.0003
			\Tstrut\Bstrut\\\cline{2-8}
			&post-processed IGA-C &3.9838&5.4367&0.2483&1.0800&1.3991&0.0051\Tstrut\Bstrut\\\thickhline
			\multirow{3}{*}{50}&Analytical& 5.0265&6.9035&0.2483&-&-&-
			\Tstrut\Bstrut\\\cline{2-8}
			&post-processed IGA-G &4.9934&6.8101&0.2483&0.6588&1.3529&0.0004
			\Tstrut\Bstrut\\\cline{2-8}
			&post-processed IGA-C &4.9797&6.7958&0.2483&0.9308&1.5591&0.0051\Tstrut\Bstrut\\
			\thickhline
		\end{tabular}
	\end{adjustbox}\label{tab:error11layers_pnt3}
\end{table}
\begin{table}[!htbp]\caption{Simply supported composite plate under a sinusoidal load with 11 layers. Out-of-plane stress state difference with respect to Pagano's solution~\cite{Pagano1970}. We compare, at $\boldsymbol{x}=(L/4,L/4,h/4)$, post-processed isogeometric collocation approach (IGA-C) and post-processed isogeometric Galerkin method (IGA-G) for a degree of approximation $p=q=6$ and 7x7 control points.}
	\vspace{0.5cm}
	\centering	
	\begin{adjustbox}{max width=\textwidth}
		\begin{tabular}{? c ? c ? c | c | c ? c | c | c ?}
			\thickhline
			\multirow{2}{*}{$\boldsymbol{S}$}&\multirow{2}{*}{\textbf{Method}}&$\sigma_{13}$ $(L/4,L/4,h/4)$
			&$\sigma_{23}$ $(L/4,L/4,h/4)$
			&$\sigma_{33}$ $(L/4,L/4,h/4)$& $\Delta(\sigma_{13})$&$\Delta(\sigma_{23})$&$\Delta(\sigma_{33})$ \Tstrut\Bstrut\\
			&&[-]&[-]&[-]&[\%]&[\%]&[\%]\Tstrut\Bstrut\\\thickhline
			\multirow{3}{*}{20}&Analytical&1.3763&2.2187&0.4209&-&-&-
			\Tstrut\Bstrut\\\cline{2-8}
			&post-processed IGA-G &1.3415&2.2104&0.4213&2.5348&0.3720&0.0931
			\Tstrut\Bstrut\\\cline{2-8}
			&post-processed IGA-C &1.3375&2.2060&0.4202&2.8209&0.5713&0.1711\Tstrut\Bstrut\\\thickhline
			\multirow{3}{*}{30}&Analytical& 2.0409&3.3489&0.4211&-&-&-
			\Tstrut\Bstrut\\\cline{2-8}
			&post-processed IGA-G &2.0122&3.3156&0.4213&1.4056&0.9939&0.0485
			\Tstrut\Bstrut\\\cline{2-8}
			&post-processed IGA-C &2.0063&3.3090&0.4202&1.6951&1.1921&0.2156\Tstrut\Bstrut\\\thickhline
			\multirow{3}{*}{40}&Analytical& 2.7094&4.4758&0.4212&-&-&-
			\Tstrut\Bstrut\\\cline{2-8}
			&post-processed IGA-G &2.6829&4.4208&0.4213&0.9777&1.2278&0.0334
			\Tstrut\Bstrut\\\cline{2-8}
			&post-processed IGA-C &2.6750&4.4120&0.4202&1.2684&1.4255&0.2306\Tstrut\Bstrut\\\thickhline
			\multirow{3}{*}{50}&Analytical& 3.3798&5.6010&0.4212&-&-&-
			\Tstrut\Bstrut\\\cline{2-8}
			&post-processed IGA-G &3.3536&5.5260&0.4213&0.7733&1.3392&0.0265
			\Tstrut\Bstrut\\\cline{2-8}
			&post-processed IGA-C &3.3438&5.5150&0.4202&1.0646&1.5366&0.2375\Tstrut\Bstrut\\
			\thickhline
		\end{tabular}
	\end{adjustbox}\label{tab:error11layers_pnt4}
\end{table}
\begin{table}[!htbp]\caption{Simply supported composite plate under a sinusoidal load with 11 layers. Out-of-plane stress state difference with respect to Pagano's solution~\cite{Pagano1970}. We compare, at $\boldsymbol{x}=(L/2,0,0)$, post-processed isogeometric collocation approach (IGA-C) and post-processed isogeometric Galerkin method (IGA-G) for a degree of approximation $p=q=6$ and 7x7 control points.}
	\vspace{0.5cm}
	\centering	
	\begin{adjustbox}{max width=\textwidth}
		\begin{tabular}{? c ? c ? c | c | c ? c | c | c ?}
			\thickhline
			\multirow{2}{*}{$\boldsymbol{S}$}&\multirow{2}{*}{\textbf{Method}}&$\sigma_{13}(L/2,0,0)$
			&$\sigma_{23}(L/2,0,0)$
			&$\sigma_{33}(L/2,0,0)$& $\Delta(\sigma_{13})$&$\Delta(\sigma_{23})$&$\Delta(\sigma_{33})$ \Tstrut\Bstrut\\
			&&[-]&[-]&[-]&[\%]&[\%]&[\%]\Tstrut\Bstrut\\\thickhline
			\multirow{3}{*}{20}&Analytical& 0.0000&5.4440&0.0000&-&-&-
			\Tstrut\Bstrut\\\cline{2-8}
			&post-processed IGA-G &0.0000&5.3558&0.0001&0.0000*&1.6191&0.0057*
			\Tstrut\Bstrut\\\cline{2-8}
			&post-processed IGA-C &0.0000&5.1533&0.0001&0.0000*&5.3394&0.0106* \Tstrut\Bstrut\\\thickhline
			\multirow{3}{*}{30}&Analytical& 0.0000& 8.2424&0.0000&-&-&-
			\Tstrut\Bstrut\\\cline{2-8}
			&post-processed IGA-G &0.0000&8.0337&0.0001&0.0000*&2.5310&0.0057*
			\Tstrut\Bstrut\\\cline{2-8}
			&post-processed IGA-C &0.0000&7.7299&0.0001&0.0000*& 6.2169& 0.0106*\Tstrut\Bstrut\\\thickhline
			\multirow{3}{*}{40}&Analytical& 0.0000& 11.0276&0.0000&-&-&-
			\Tstrut\Bstrut\\\cline{2-8}
			&post-processed IGA-G &0.0000&10.7117&0.0001&0.0000*&2.8653&0.0057*
			\Tstrut\Bstrut\\\cline{2-8}
			&post-processed IGA-C &0.0000&10.3066&0.0001&0.0000*&6.5386& 0.0106*\Tstrut\Bstrut\\\thickhline
			\multirow{3}{*}{50}&Analytical& 0.0000&13.8069&0.0000&-&-&-
			\Tstrut\Bstrut\\\cline{2-8}
			&post-processed IGA-G &0.0000&13.3896&0.0001&0.0000*&3.0229&0.0057*
			\Tstrut\Bstrut\\\cline{2-8}
			&post-processed IGA-C &0.0000&12.8832&0.0001&0.0000*&6.6902& 0.0106*\Tstrut\Bstrut\\
			\thickhline
		\end{tabular}
	\end{adjustbox}\label{tab:error11layers_pnt5}
\end{table}
\begin{table}[!htbp]\caption{Simply supported composite plate under a sinusoidal load with 11 layers. Out-of-plane stress state difference with respect to Pagano's solution~\cite{Pagano1970}. We compare, at $\boldsymbol{x}=(L/2,0,h/4)$, post-processed isogeometric collocation approach (IGA-C) and post-processed isogeometric Galerkin method (IGA-G) for a degree of approximation $p=q=6$ and 7x7 control points.}
	\vspace{0.5cm}
	\centering	
	\begin{adjustbox}{max width=\textwidth}
		\begin{tabular}{? c ? c ? c | c | c ? c | c | c ?}
			\thickhline
			\multirow{2}{*}{$\boldsymbol{S}$}&\multirow{2}{*}{\textbf{Method}}&$\sigma_{13}(L/2,0,h/4)$
			&$\sigma_{23}(L/2,0,h/4)$
			&$\sigma_{33}(L/2,0,h/4)$& $\Delta(\sigma_{13})$&$\Delta(\sigma_{23})$&$\Delta(\sigma_{33})$ \Tstrut\Bstrut\\
			&&[-]&[-]&[-]&[\%]&[\%]&[\%]\Tstrut\Bstrut\\\thickhline
			\multirow{3}{*}{20}&Analytical& 0.0000&4.4373&0.0000&-&-&-
			\Tstrut\Bstrut\\\cline{2-8}
			&post-processed IGA-G &0.0000&4.3456&-0.0047&0.0000*&2.0680&0.4720*
			\Tstrut\Bstrut\\\cline{2-8}
			&post-processed IGA-C &0.0000&4.1806&-0.0087&0.0000*&5.7856&0.8691* \Tstrut\Bstrut\\\thickhline
			\multirow{3}{*}{30}&Analytical& 0.0000& 6.6978&0.0000&-&-&-
			\Tstrut\Bstrut\\\cline{2-8}
			&post-processed IGA-G &0.0000&6.5183&-0.0047&0.0000*&2.6794&0.4720*
			\Tstrut\Bstrut\\\cline{2-8}
			&post-processed IGA-C &0.0000&6.2709&-0.0087&0.0000*& 6.3738& 0.8691*\Tstrut\Bstrut\\\thickhline
			\multirow{3}{*}{40}&Analytical& 0.0000& 8.9516&0.0000&-&-&-
			\Tstrut\Bstrut\\\cline{2-8}
			&post-processed IGA-G &0.0000&8.6911&-0.0047&0.0000*&2.9094&0.4720*
			\Tstrut\Bstrut\\\cline{2-8}
			&post-processed IGA-C &0.0000&8.3612&-0.0087&0.0000*&6.5950& 0.8691*\Tstrut\Bstrut\\\thickhline
			\multirow{3}{*}{50}&Analytical& 0.0000&11.2021&0.0000&-&-&-
			\Tstrut\Bstrut\\\cline{2-8}
			&post-processed IGA-G &0.0000&10.8639&-0.0047&0.0000*&3.0188&0.4720*
			\Tstrut\Bstrut\\\cline{2-8}
			&post-processed IGA-C &0.0000&10.4515&-0.0087&0.0000*&6.7003& 0.8691*\Tstrut\Bstrut\\
			\thickhline
		\end{tabular}
	\end{adjustbox}\label{tab:error11layers_pnt6}
\end{table}
\begin{table}[!htbp]\caption{Simply supported composite plate under a sinusoidal load with 34 layers. Out-of-plane stress state difference with respect to Pagano's solution~\cite{Pagano1970}. We compare, at $\boldsymbol{x}=(0,L/2,0)$, post-processed isogeometric collocation approach (IGA-C) and post-processed isogeometric Galerkin method (IGA-G) for a degree of approximation $p=q=6$ and 7x7 control points.}
	\vspace{0.5cm}
	\centering	
	\begin{adjustbox}{max width=\textwidth}
		\begin{tabular}{? c ? c ? c | c | c ? c | c | c ?}
			\thickhline
			\multirow{2}{*}{$\boldsymbol{S}$}&\multirow{2}{*}{\textbf{Method}}&$\sigma_{13}(0,L/2,0)$
			&$\sigma_{23}(0,L/2,0)$
			&$\sigma_{33}(0,L/2,0)$& $\Delta(\sigma_{13})$&$\Delta(\sigma_{23})$&$\Delta(\sigma_{33})$ \Tstrut\Bstrut\\
			&&[-]&[-]&[-]&[\%]&[\%]&[\%]\Tstrut\Bstrut\\\thickhline
			\multirow{3}{*}{20}&Analytical& 4.7476&0.0000&0.0000&-&-&-
			\Tstrut\Bstrut\\\cline{2-8}
			&post-processed IGA-G &4.6422&0.0000&0.0000&2.2199&0.0000*&0.0019*
			\Tstrut\Bstrut\\\cline{2-8}
			&post-processed IGA-C &4.4689&0.0000&0.0000&5.8699&0.0000*&0.0034*\Tstrut\Bstrut\\\thickhline
			\multirow{3}{*}{30}&Analytical& 7.1411&0.0000&0.0000&-&-&-
			\Tstrut\Bstrut\\\cline{2-8}
			&post-processed IGA-G &6.9633&0.0000&0.0000&2.4890&0.0000*&0.0019*
			\Tstrut\Bstrut\\\cline{2-8}
			&post-processed IGA-C &6.7034&0.0000&0.0000&6.1289&0.0000*&0.0034*\Tstrut\Bstrut\\\thickhline
			\multirow{3}{*}{40}&Analytical& 9.5307&0.0000&0.0000&-&-&-
			\Tstrut\Bstrut\\\cline{2-8}
			&post-processed IGA-G &9.2844&0.0000&0.0000&2.5839&0.0000*&0.0019*
			\Tstrut\Bstrut\\\cline{2-8}
			&post-processed IGA-C &8.9378&0.0000&0.0000&6.2203&0.0000*&0.0034*\Tstrut\Bstrut\\\thickhline
			\multirow{3}{*}{50}&Analytical& 11.9187&0.0000&0.0000&-&-&-
			\Tstrut\Bstrut\\\cline{2-8}
			&post-processed IGA-G &11.6055&0.0000&0.0000&2.6280&0.0000*&0.0019*
			\Tstrut\Bstrut\\\cline{2-8}
			&post-processed IGA-C &11.1723&0.0000&0.0000&6.2627&0.0000*&0.0034*\Tstrut\Bstrut\\
			\thickhline
		\end{tabular}
	\end{adjustbox}\label{tab:error34layers_pnt1}
\end{table}
\begin{table}[!htbp]\caption{Simply supported composite plate under a sinusoidal load with 34 layers. Out-of-plane stress state difference with respect to Pagano's solution~\cite{Pagano1970}. We compare, at $\boldsymbol{x}=(0,L/2,h/4)$, post-processed isogeometric collocation approach (IGA-C) and post-processed isogeometric Galerkin method (IGA-G) for a degree of approximation $p=q=6$ and 7x7 control points.}
	\vspace{0.5cm}
	\centering	
	\begin{adjustbox}{max width=\textwidth}
		\begin{tabular}{? c ? c ? c | c | c ? c | c | c ?}
			\thickhline
			\multirow{2}{*}{$\boldsymbol{S}$}&\multirow{2}{*}{\textbf{Method}}&$\sigma_{13}(0,L/2,h/4)$
			&$\sigma_{23}(0,L/2,h/4)$
			&$\sigma_{33}(0,L/2,h/4)$& $\Delta(\sigma_{13})$&$\Delta(\sigma_{23})$&$\Delta(\sigma_{33})$ \Tstrut\Bstrut\\
			&&[-]&[-]&[-]&[\%]&[\%]&[\%]\Tstrut\Bstrut\\\thickhline
			\multirow{3}{*}{20}&Analytical& 3.7058&0.0000&0.0000&-&-&-
			\Tstrut\Bstrut\\\cline{2-8}
			&post-processed IGA-G &3.5969&0.0000&-0.0046&2.9384&0.0000*&0.4610*
			\Tstrut\Bstrut\\\cline{2-8}
			&post-processed IGA-C &3.4624&0.0000&-0.0085&6.5684&0.0000*&0.8461*\Tstrut\Bstrut\\\thickhline
			\multirow{3}{*}{30}&Analytical& 5.5563&0.0000&0.0000&-&-&-
			\Tstrut\Bstrut\\\cline{2-8}
			&post-processed IGA-G &5.3954&0.0000&-0.0046&2.8963&0.0000*&0.4610*
			\Tstrut\Bstrut\\\cline{2-8}
			&post-processed IGA-C &5.1936&0.0000&-0.0085&6.5279&0.0000*&0.8461*\Tstrut\Bstrut\\\thickhline
			\multirow{3}{*}{40}&Analytical& 7.4073&0.0000&0.0000&-&-&-
			\Tstrut\Bstrut\\\cline{2-8}
			&post-processed IGA-G &7.1938&0.0000&-0.0046&2.8821&0.0000*&0.4610*
			\Tstrut\Bstrut\\\cline{2-8}
			&post-processed IGA-C &6.9248&0.0000&-0.0085&6.5143&0.0000*&0.8461*\Tstrut\Bstrut\\\thickhline
			\multirow{3}{*}{50}&Analytical& 9.2585&0.0000&0.0000&-&-&-
			\Tstrut\Bstrut\\\cline{2-8}
			&post-processed IGA-G &8.9923&0.0000&-0.0046&2.8757&0.0000*&0.4610*
			\Tstrut\Bstrut\\\cline{2-8}
			&post-processed IGA-C &8.6560&0.0000&-0.0085&6.5081&0.0000*&0.8461*\Tstrut\Bstrut\\
			\thickhline
		\end{tabular}
	\end{adjustbox}\label{tab:error34layers_pnt2}
\end{table}
\begin{table}[!htbp]\caption{Simply supported composite plate under a sinusoidal load with 34 layers. Out-of-plane stress state difference with respect to Pagano's solution~\cite{Pagano1970}. We compare, at $\boldsymbol{x}=(L/4,L/4,0)$, post-processed isogeometric collocation approach (IGA-C) and post-processed isogeometric Galerkin method (IGA-G) for a degree of approximation $p=q=6$ and 7x7 control points.}
	\vspace{0.5cm}
	\centering	
	\begin{adjustbox}{max width=\textwidth}
		\begin{tabular}{? c ? c ? c | c | c ? c | c | c ?}
			\thickhline
			\multirow{2}{*}{$\boldsymbol{S}$}&\multirow{2}{*}{\textbf{Method}}&$\sigma_{13}$ $(L/4,L/4,0)$
			&$\sigma_{23}$ $(L/4,L/4,0)$
			&$\sigma_{33}$ $(L/4,L/4,0)$& $\Delta(\sigma_{13})$&$\Delta(\sigma_{23})$&$\Delta(\sigma_{33})$ \Tstrut\Bstrut\\
			&&[-]&[-]&[-]&[\%]&[\%]&[\%]\Tstrut\Bstrut\\\thickhline
			\multirow{3}{*}{20}&Analytical&2.3738&2.3746&0.2494&-&-&-
			\Tstrut\Bstrut\\\cline{2-8}
			&post-processed IGA-G &2.3606&2.3609&0.2494&0.5570&0.5767&0.0029
			\Tstrut\Bstrut\\\cline{2-8}
			&post-processed IGA-C &2.3550&2.3553&0.2495&0.7930&0.8127&0.0046\Tstrut\Bstrut\\\thickhline
			\multirow{3}{*}{30}&Analytical& 3.5705&3.5713&0.2494&-&-&-
			\Tstrut\Bstrut\\\cline{2-8}
			&post-processed IGA-G &3.5409&3.5413&0.2494&0.8307&0.8404&0.0001
			\Tstrut\Bstrut\\\cline{2-8}
			&post-processed IGA-C &3.5325&3.5329&0.2495&1.0660&1.0757&0.0019\Tstrut\Bstrut\\\thickhline
			\multirow{3}{*}{40}&Analytical& 4.7653&4.7663&0.2494&-&-&-
			\Tstrut\Bstrut\\\cline{2-8}
			&post-processed IGA-G &4.7212&4.7218&0.2494&0.9272&0.9332&0.0002
			\Tstrut\Bstrut\\\cline{2-8}
			&post-processed IGA-C &4.7099&4.7106&0.2495&1.1623&1.1683&0.0015\Tstrut\Bstrut\\\thickhline
			\multirow{3}{*}{50}&Analytical& 5.9594&5.9604&0.2494&-&-&-
			\Tstrut\Bstrut\\\cline{2-8}
			&post-processed IGA-G &5.9014&5.9022&0.2494&0.9720&0.9763&0.0003
			\Tstrut\Bstrut\\\cline{2-8}
			&post-processed IGA-C &5.8874&5.8882&0.2495&1.2071&1.2113&0.0014\Tstrut\Bstrut\\
			\thickhline
		\end{tabular}
	\end{adjustbox}\label{tab:error34layers_pnt3}
\end{table}
\begin{table}[!htbp]\caption{Simply supported composite plate under a sinusoidal load with 34 layers. Out-of-plane stress state difference with respect to Pagano's solution~\cite{Pagano1970}. We compare, at $\boldsymbol{x}=(L/4,L/4,h/4)$, post-processed isogeometric collocation approach (IGA-C) and post-processed isogeometric Galerkin method (IGA-G) for a degree of approximation $p=q=6$ and 7x7 control points.}
	\vspace{0.5cm}
	\centering	
	\begin{adjustbox}{max width=\textwidth}
		\begin{tabular}{? c ? c ? c | c | c ? c | c | c ?}
			\thickhline
			\multirow{2}{*}{$\boldsymbol{S}$}&\multirow{2}{*}{\textbf{Method}}&$\sigma_{13}$ $(L/4,L/4,h/4)$
			&$\sigma_{23}$ $(L/4,L/4,h/4)$
			&$\sigma_{33}$ $(L/4,L/4,h/4)$& $\Delta(\sigma_{13})$&$\Delta(\sigma_{23})$&$\Delta(\sigma_{33})$ \Tstrut\Bstrut\\
			&&[-]&[-]&[-]&[\%]&[\%]&[\%]\Tstrut\Bstrut\\\thickhline
			\multirow{3}{*}{20}&Analytical&1.8529&1.7370&0.4212&-&-&-
			\Tstrut\Bstrut\\\cline{2-8}
			&post-processed IGA-G &1.8291&1.7154&0.4217&1.2830&1.2387&0.1167
			\Tstrut\Bstrut\\\cline{2-8}
			&post-processed IGA-C &1.8249&1.7113&0.4206&1.5130&1.4777&0.1476\Tstrut\Bstrut\\\thickhline
			\multirow{3}{*}{30}&Analytical& 2.7782&2.6027&0.4215&-&-&-
			\Tstrut\Bstrut\\\cline{2-8}
			&post-processed IGA-G &2.7437&2.5732&0.4217&1.2402&1.1343&0.0481
			\Tstrut\Bstrut\\\cline{2-8}
			&post-processed IGA-C &2.7373&2.5669&0.4206&1.4703&1.3735&0.2160\Tstrut\Bstrut\\\thickhline
			\multirow{3}{*}{40}&Analytical& 3.7037&3.4690&0.4216&-&-&-
			\Tstrut\Bstrut\\\cline{2-8}
			&post-processed IGA-G &3.6583&3.4309&0.4217&1.2258&1.0978&0.0241
			\Tstrut\Bstrut\\\cline{2-8}
			&post-processed IGA-C &3.6497&3.4226&0.4206&1.4559&1.3371&0.2399\Tstrut\Bstrut\\\thickhline
			\multirow{3}{*}{50}&Analytical& 4.6293&4.3355&0.4217&-&-&-
			\Tstrut\Bstrut\\\cline{2-8}
			&post-processed IGA-G &4.5728&4.2886&0.4217&1.2192&1.0809&0.0130
			\Tstrut\Bstrut\\\cline{2-8}
			&post-processed IGA-C &4.5622&4.2782&0.4206&1.4494&1.3203&0.2510\Tstrut\Bstrut\\
			\thickhline
		\end{tabular}
	\end{adjustbox}\label{tab:error34layers_pnt4}
\end{table}
\begin{table}[!htbp]\caption{Simply supported composite plate under a sinusoidal load with 34 layers. Out-of-plane stress state difference with respect to Pagano's solution~\cite{Pagano1970}. We compare, at $\boldsymbol{x}=(L/2,0,0)$, post-processed isogeometric collocation approach (IGA-C) and post-processed isogeometric Galerkin method (IGA-G) for a degree of approximation $p=q=6$ and 7x7 control points.}
	\vspace{0.5cm}
	\centering	
	\begin{adjustbox}{max width=\textwidth}
		\begin{tabular}{? c ? c ? c | c | c ? c | c | c ?}
			\thickhline
			\multirow{2}{*}{$\boldsymbol{S}$}&\multirow{2}{*}{\textbf{Method}}&$\sigma_{13}(L/2,0,0)$
			&$\sigma_{23}(L/2,0,0)$
			&$\sigma_{33}(L/2,0,0)$& $\Delta(\sigma_{13})$&$\Delta(\sigma_{23})$&$\Delta(\sigma_{33})$ \Tstrut\Bstrut\\
			&&[-]&[-]&[-]&[\%]&[\%]&[\%]\Tstrut\Bstrut\\\thickhline
			\multirow{3}{*}{20}&Analytical& 0.0000&4.7492&0.0000&-&-&-
			\Tstrut\Bstrut\\\cline{2-8}
			&post-processed IGA-G &0.0000&4.6428&0.0000&0.0000*&2.2394&0.0018*
			\Tstrut\Bstrut\\\cline{2-8}
			&post-processed IGA-C &0.0000&4.4695&0.0000&0.0000*&5.8886&0.0032* \Tstrut\Bstrut\\\thickhline
			\multirow{3}{*}{30}&Analytical& 0.0000&7.1427&0.0000&-&-&-
			\Tstrut\Bstrut\\\cline{2-8}
			&post-processed IGA-G &0.0000&6.9642&0.0000&0.0000*&2.4986&0.0018*
			\Tstrut\Bstrut\\\cline{2-8}
			&post-processed IGA-C &0.0000&6.7043&0.0000&0.0000*&6.1382&0.0032*\Tstrut\Bstrut\\\thickhline
			\multirow{3}{*}{40}&Analytical& 0.0000&9.5325&0.0000&-&-&-
			\Tstrut\Bstrut\\\cline{2-8}
			&post-processed IGA-G &0.0000&9.2856&0.0000&0.0000*&2.5899&0.0018*
			\Tstrut\Bstrut\\\cline{2-8}
			&post-processed IGA-C &0.0000&8.9390&0.0000&0.0000*&6.2260&0.0032*\Tstrut\Bstrut\\\thickhline
			\multirow{3}{*}{50}&Analytical& 0.0000&11.9208&0.0000&-&-&-
			\Tstrut\Bstrut\\\cline{2-8}
			&post-processed IGA-G &0.0000&11.6070&0.0000&0.0000*&2.6322&0.0018*
			\Tstrut\Bstrut\\\cline{2-8}
			&post-processed IGA-C &0.0000&11.1738&0.0000&0.0000*&6.2668&0.0032*\Tstrut\Bstrut\\
			\thickhline
		\end{tabular}
	\end{adjustbox}\label{tab:error34layers_pnt5}
\end{table}
\begin{table}[!htbp]\caption{Simply supported composite plate under a sinusoidal load with 34 layers. Out-of-plane stress state difference with respect to Pagano's solution~\cite{Pagano1970}. We compare, at $\boldsymbol{x}=(L/2,0,h/4)$, post-processed isogeometric collocation approach (IGA-C) and post-processed isogeometric Galerkin method (IGA-G) for a degree of approximation $p=q=6$ and 7x7 control points.}
	\vspace{0.5cm}
	\centering	
	\begin{adjustbox}{max width=\textwidth}
		\begin{tabular}{? c ? c ? c | c | c ? c | c | c ?}
			\thickhline
			\multirow{2}{*}{$\boldsymbol{S}$}&\multirow{2}{*}{\textbf{Method}}&$\sigma_{13}(L/2,0,h/4)$
			&$\sigma_{23}(L/2,0,h/4)$
			&$\sigma_{33}(L/2,0,h/4)$& $\Delta(\sigma_{13})$&$\Delta(\sigma_{23})$&$\Delta(\sigma_{33})$ \Tstrut\Bstrut\\
			&&[-]&[-]&[-]&[\%]&[\%]&[\%]\Tstrut\Bstrut\\\thickhline
			\multirow{3}{*}{20}&Analytical& 0.0000&3.4739&0.0000&-&-&-
			\Tstrut\Bstrut\\\cline{2-8}
			&post-processed IGA-G &0.0000&3.3737&-0.0046&0.0000*&2.8853&0.4583*
			\Tstrut\Bstrut\\\cline{2-8}
			&post-processed IGA-C &0.0000&3.2480&-0.0084&0.0000*&6.5031&0.8413* \Tstrut\Bstrut\\\thickhline
			\multirow{3}{*}{30}&Analytical& 0.0000&5.2054&0.0000&-&-&-
			\Tstrut\Bstrut\\\cline{2-8}
			&post-processed IGA-G &0.0000&5.0605&-0.0046&0.0000*&2.7826&0.4583*
			\Tstrut\Bstrut\\\cline{2-8}
			&post-processed IGA-C &0.0000&4.8720&-0.0084&0.0000*&6.4042&0.8413*\Tstrut\Bstrut\\\thickhline
			\multirow{3}{*}{40}&Analytical& 0.0000&6.9380&0.0000&-&-&-
			\Tstrut\Bstrut\\\cline{2-8}
			&post-processed IGA-G &0.0000&6.7474&-0.0046&0.0000*&2.7468&0.4583*
			\Tstrut\Bstrut\\\cline{2-8}
			&post-processed IGA-C &0.0000&6.4960&-0.0084&0.0000*&6.3697&0.8413*\Tstrut\Bstrut\\\thickhline
			\multirow{3}{*}{50}&Analytical& 0.0000&8.6710&0.0000&-&-&-
			\Tstrut\Bstrut\\\cline{2-8}
			&post-processed IGA-G &0.0000&8.4342&-0.0046&0.0000*&2.7302&0.4583*
			\Tstrut\Bstrut\\\cline{2-8}
			&post-processed IGA-C &0.0000&8.1200&-0.0084&0.0000*&6.3537&0.8413*\Tstrut\Bstrut\\
			\thickhline
		\end{tabular}
	\end{adjustbox}\label{tab:error34layers_pnt6}
\end{table}

\newpage
\section{Conclusions}\label{sec:conclusions}

Moving from the equilibrium-based post-processing technique that we have recently proposed in the context of 3D solid plates approximated by isogeometric Galerkin \cite{Dufour2018} or collocation \cite{Patton2019} methods, in this paper, we have considered the application of such an approach for the accurate and inexpensive recovery of interlaminar stresses in IGA Kirchhoff plates. The adopted method has been shown to be particularly effective in this framework for both Galerkin and collocation methods, also in the case of an even number of layers and, therefore, of non-symmetric ply distributions.
The fundamental ingredients to obtain such good results are the high accuracy and regularity granted by high-order IGA methods even with coarse meshes, and our numerical tests for a simple gometry like the one involved in Pagano test case have shown that even a mesh constituted by a single sixth-order element is able to provide very good results in terms of in-plane and out-of-plane stresses, for both Galerkin and collocation discretizations.
Extensive numerical experiments have confirmed the high efficiency of the proposed approach.

Finally, the extension to curved geometries is currently under investigation and very convincing preliminary results have already been obtained, while the application to nonlinear problems will be the subject of future research.

\section*{Acknowledgements}
This work was partially supported by Ministero dell'Istruzione, dell'Universit\`a e della Ricerca through the project ``XFAST-SIMS: Extra fast and accurate simulation of complex structural systems'', within the program Progetti di ricerca di Rilevante Interesse Nazionale (PRIN). P. Antol\'in was also partially supported by the European Research Council through the H2020 ERC Advanced Grant 2015 n. 694515 CHANGE, and by the Swiss National Science Foundation through the project ``Design-through-Analysis (of PDEs): the litmus test'' n. 40B2-0 187094 (BRIDGE Discovery 2019). J. Kiendl was partially supported by the European Research Council through the H2020 ERC Consolidator Grant 2019 n. 864482 FDM$^2$.

% ...

\newpage
\bibliographystyle{unsrt}
\bibliography{references}

\end{document}